\documentclass[10pt,twocolumn,twoside]{IEEEtran}

\usepackage[T1]{fontenc}
\usepackage{units}
\usepackage{amsthm}
\usepackage{amsmath}
\usepackage{amssymb}
\usepackage{enumerate,color,cite}
\usepackage{ragged2e}
\usepackage{graphicx}
\usepackage{algorithm, algorithmic}
\usepackage[hidelinks]{hyperref}

\graphicspath{{./Figs/}}

\makeatletter

\theoremstyle{plain}
\newtheorem{thm}{\theoremname}
\theoremstyle{plain}
\newtheorem{lem}{\lemmaname}
\theoremstyle{plain}

\theoremstyle{definition}
\newtheorem{defn}{\definitionname}
\theoremstyle{remark}
\newtheorem{rem}{\remarkname}

\newtheorem{assumption}{Assumption}[section]

\makeatother

\providecommand{\corollaryname}{Corollary}
\providecommand{\definitionname}{Definition}
\providecommand{\lemmaname}{Lemma}
\providecommand{\remarkname}{Remark}
\providecommand{\theoremname}{Theorem}

\usepackage{caption}
\newcommand*{\TitleFont}{%
      \fontsize{22}{26}%
      \selectfont}

\newcommand{\K}{\mathcal{K}}
\newcommand{\Kinf}{\K_{\infty}}
\newcommand{\LL}{\mathcal{L}}
\newcommand{\KL}{\mathcal{KL}}
\newcommand{\KdL}{\mathcal{K \cdot L}}

\begin{document}

\title{\TitleFont Robust Stability of Optimization-based State Estimation}


\author{Wuhua~Hu%
\thanks{W. Hu is with the Institute for Infocomm Research, Agency for Science,
Technology and Research (A$^\star$STAR), Singapore,  E-mail: huwh@i2r.a-star.edu.sg}}

\maketitle
\begin{abstract}
Optimization-based state estimation is useful for nonlinear or constrained
dynamic systems for which few general methods with established properties
are available. The two fundamental forms are moving horizon estimation
(MHE) which uses the nearest measurements within a moving time horizon,
and its theoretical ideal, full information estimation (FIE) which
uses all measurements up to the time of estimation. Despite extensive
studies, the stability analyses of FIE and MHE for discrete-time nonlinear
systems with bounded process and measurement disturbances, remain
an open challenge. This work aims to provide a systematic solution
for the challenge. First, we prove that FIE is robustly globally asymptotically
stable (RGAS) if the cost function admits a property mimicking the
incremental input/output-to-state stability (i-IOSS) of the system
and has a sufficient sensitivity to the uncertainty in the initial
state. Second, we establish an explicit link from the RGAS of FIE
to that of MHE, and use it to show that MHE is RGAS under enhanced
conditions if the moving horizon is long enough to suppress the propagation
of uncertainties. The theoretical results imply flexible MHE designs
with assured robust stability for a broad class of i-IOSS systems.
Numerical experiments on linear and nonlinear systems are used to
illustrate the designs and support the findings.
\end{abstract}

\thispagestyle{empty}

\section{Introduction}

Optimization-based state estimation refers to an estimation method
that estimates the state of a system via an optimization approach,
in which the optimization utilizes all or a subset of the information
available about the system up to the time of estimation. It has advantages
in handling nonlinear or constrained systems for which few general
state estimation methods with established properties are available
\cite{rawlings2014moving}. Full information estimation (FIE) is an
ideal form of optimization-based estimation which uses all measurements
up to the time of estimation. In the absence of constraints, FIE is
equivalent to Kalman filtering (KF) when the system is linear time-invariant
and the cost function has an appropriate quadratic form \cite{rawlings2014moving}.
Since the measurements increase with time, it is impractical to implement
FIE. This motivates the development of moving horizon estimation (MHE)
as its practical approximation which uses only the latest batch of
measurements to do the estimation. The idea of MHE dates back to 1960's
\cite{jazwinski1968limited} which was motivated for making KF robust
to modeling errors. However, it is not until recently that the idea
is gradually developed into a field, i.e., the field of MHE \cite{muske1995nonlinear,allgower1999nonlinear,rawlings2012optimization,rawlings2014moving}.
The recent developments include MHE theoretical and applied researches
which investigate the stability and the implementation issues of MHE.

Theoretical research has been concentrated on the stability conditions
of MHE. Early research assumed linear systems \cite{muske1993receding,rao2001constrained,alessandri2003receding},
and later nonlinear systems \cite{michalska1995moving,alessandri1999neural,rao2003constrained,alessandri2008moving}.
{
Part of the stability results were obtained by assuming the presence
of measurement disturbances but the absence of process disturbances,
e.g., \cite{alessandri1999neural}, and most were obtained by assuming
the presence of both disturbances, e.g., \cite{muske1993receding,michalska1995moving,rao2001constrained,rao2003constrained,alessandri2003receding,alessandri2008moving}.
}%
And the stability results were derived based on different formulations
of the problems. For example, references \cite{muske1993receding,rao2001constrained,alessandri2003receding,alessandri2010advances,liu2013moving}
considered either linear or nonlinear systems and each assumed a quadratic
cost function that accounts for a quadratic arrival cost and quadratic
penalties on the measurement fitting errors. Whereas, references \cite{rao2003constrained,rawlings2012optimization,hu2015optimization,ji2016robust,muller2017nonlinear}
considered general nonlinear systems and assumed a general form of
cost functions that are not necessarily in a quadratic form.

The recent review made in \cite{rawlings2012optimization} provides
a concise and general view of the problem which relies on the concept
of incremental input/output-to-state stability (i-IOSS; refer to Definition \ref{def: i-IOSS} in Section \ref{sec: Notation-and-Preliminaries}) for detectability
of nonlinear systems \cite{sontag1997output} and the concept of robust
global asymptotic stability (RGAS) for robust stability of a state
estimator \cite{rawlings2009model,rawlings2012optimization}. The
review revealed two major challenges that were open in the field \cite{rawlings2012optimization,rawlings2014moving}:
(i) the search for conditions and a proof of the RGAS of MHE in the
presence of bounded disturbances, and (ii) the development of suboptimal
MHE that enables an efficient computation of the solution. As an initial
step to tackling challenge (i), reference \cite{hu2015optimization}
identified a broad class of cost functions that ensure the RGAS of
FIE, and the cost functions were shown to admit a more specific form
for a class of i-IOSS systems considered in \cite{ji2016robust}.
{
The implication to the RGAS of MHE was further investigated in \cite{muller2017nonlinear}
based on the results of \cite{ji2016robust}, which showed that MHE
is RGAS if the same conditions are enhanced properly. Moreover, in \cite{muller2017nonlinear} the convergence of MHE for convergent disturbances was proved under the enhanced conditions, and the MHE was shown to be RGAS even if the cost function does not have max terms which are needed in the stability analyses of FIE \cite{ji2016robust}.
}

Other relevant progress was reported in \cite{liu2013moving}, which
assumed a quadratic cost function and used a nonlinear deterministic
observer to generate useful constraints so that the MHE results in
bounded estimation errors under certain conditions. Some earlier developments
are also available in \cite{alessandri2008moving,alessandri2010advances},
which assumed a quadratic cost function for MHE and also an observability
and some Lipschitz conditions on the system. The other developments
of MHE are mainly in applied research, which mostly aimed at reducing
the online computational complexity of MHE for applications in large
dimensional and nonlinear systems \cite{rawlings2014moving}. Interested
readers are referred to \cite{zavala2008fast,kuhl2011real,lopez2012moving,wynn2014convergence,alessandri2016moving_cdc}
and the references therein for relevant examples.

This paper follows the general view of MHE developed in \cite{rawlings2012optimization},
and aims to present a systematic solution for the aforementioned open
challenge (i). It significantly extends our conference paper \cite{hu2015optimization},
which identified sufficient conditions for the RGAS of FIE. The major
new contributions are the following. First, an explicit link is established
from the RGAS of FIE to that of MHE. Second, based on the link we
prove the RGAS of MHE under both general and specialized conditions,
depending on the stability property of a system. The convergence of MHE to the true state is also established for a system subjected to disturbances that converge to zero. Third, two numerical
examples are developed with rigorous analyses to support the theoretical
findings. The numerical results are compared with those obtained by
KF in the linear case and extended Kalman filter (EKF) \cite{ljung1979asymptotic}
in the nonlinear case, demonstrating the advantages of MHE in the
considered situations. 

{
Since the discrete-time system and
the MHE are assumed to have general forms, the theoretical results
imply flexible MHE designs with assured robust stability for a broad
class of systems. This constitutes a key difference from the recent
stability results obtained in \cite{muller2017nonlinear} which are applicable to a subset of the considered
systems.
}

The rest of the paper is organized as follows. Section \ref{sec: Notation-and-Preliminaries}
introduces the notation and some preliminaries. Section \ref{sec: optimization-based estimation}
defines the ideal and the practical forms of optimization-based state
estimation, i.e., FIE and MHE. Section \ref{sec: RGAS of FIE} proves
the robust stability of FIE under general and then specialized conditions.
Section \ref{sec: RGAS of the MHE} reveals its implication to MHE
and subsequently proves the robust stability of MHE under enhanced
conditions. Convergence of the MHE is also proved for disturbances that are convergent to zero. Section \ref{sec:numerical example} presents two numerical
examples to illustrate the flexible MHE designs. Section \ref{sec:Discussion} provides a brief discussion on ways to tackle the computational challenge of MHE. Finally, Section
\ref{sec:Conclusion} concludes the paper with a remark on the future
work.

\section{Notation and Preliminaries \label{sec: Notation-and-Preliminaries}}

The notation mostly follows the convention in \cite{rawlings2012optimization}.
The symbols $\mathbb{R}$, $\mathbb{R}_{\ge0}$ and $\mathbb{I}_{\ge0}$
denote the sets of real numbers, nonnegative real numbers and nonnegative
integers, respectively, and $\mathbb{I}_{a:b}$ denotes the set of
integers from $a$ to $b$. The constraints $t\ge0$ and $t\in\mathbb{I}_{\ge0}$
are used interchangeably to refer to the set of discrete times. The
symbol $\left|\cdot\right|$ denotes the Euclidean norm of a vector
or the 2-norm of a matrix, depending on the argument. The bold symbol
$\boldsymbol{x}_{a:b}$, denotes a sequence of vector-valued variables
$\{x_{a},\,x_{a+1},\,...,\thinspace x_{b}\}$, and with a function
$f$ acting on a vector $x$, $f(\boldsymbol{x}_{a:b})$ stands for
the sequence of function values $\{f(x_{a}),\,f(x_{a+1}),\,...,\thinspace f(x_{b})\}$.
The notation $\left\Vert \boldsymbol{x}_{a:b}\right\Vert $ refers
to $\max_{a\le i\le b}\left|x_{i}\right|$ if $a\le b$ and to 0 if
$a>b$. Throughout the paper, $t$ refers to a discrete time, and
as a subscript it indicates dependence on time $t$. Whereas, the
subscripts or superscripts $x$, $w$ and $v$ are used exclusively
to indicate a function or variable that is associated with the state
($x$), disturbance ($w$) or measurement noise ($v$), and they do
not imply dependence relationships. The frequently used $\K$, $\Kinf$,
$\LL$ and $\KL$ functions are defined as follows.
\begin{defn}
($\K$, $\Kinf$, $\LL$ and $\KL$ functions) A function $\alpha:\mathbb{R}_{\ge0}\to\mathbb{R}_{\ge0}$
is a $\K$ function if it is continuous, zero at zero, and strictly
increasing, and a $\Kinf$ function if $\alpha$ is a $\K$ function
and satisfies $\alpha(s)\to\infty$ as $s\to\infty$. A function $\varphi:\mathbb{R}_{\ge0}\to\mathbb{R}_{\ge0}$
is a $\LL$ function if it is continuous, nonincreasing and satisfies
$\varphi(t)\to0$ as $t\to\infty$. A function $\beta:\mathbb{R}_{\ge0}\times\mathbb{R}_{\ge0}\to\mathbb{R}_{\ge0}$
is a $\KL$ function if, for each $t\ge0$, $\beta(\cdot,t)$ is
a $\K$ function and for each $s\ge0$, $\beta(s,\cdot)$ is a $\LL$
function. 
\end{defn}
The following properties of the $\K$ and $\KL$ functions will be
used in later analyses. 
\begin{lem}
\cite{khalil2002nonlinear,rawlings2012optimization} Given a $\K$
function $\alpha$ and a $\KL$ function $\beta$, the following
holds for all $a_{i}\in\mathbb{R}_{\ge0}$, $i\in\mathbb{I}_{1:n}$,
and all $t\in\mathbb{R}_{\ge0}$, 
\[
\alpha\left(\sum_{i=1}^{n}a_{i}\right)\le\sum_{i=1}^{n}\alpha(na_{i}),\,\,\beta\left(\sum_{i=1}^{n}a_{i},t\right)\le\sum_{i=1}^{n}\beta(na_{i},t).
\]
\end{lem}
\begin{defn}
\label{def: i-IOSS}(i-IOSS \cite{rawlings2012optimization,sontag1997output})
The system $x_{t+1}=f(x_{t},w_{t})$, $y_{t}=h(x_{t})$ is i-IOSS
if there exist functions $\beta\in\KL$ and $\alpha_{1},\alpha_{2}\in\K$
such that for every two initial states $x_{0}^{(1)},x_{0}^{(2)}$,
and two sequences of disturbances $\boldsymbol{w}_{0:t-1}^{(1)},\boldsymbol{w}_{0:t-1}^{(2)}$,
the following inequality holds for all $t\in\mathbb{I}_{\ge0}$: 
\begin{align}
 & \left|x_{t}(x_{0}^{(1)},\boldsymbol{w}_{0:t-1}^{(1)})-x_{t}(x_{0}^{(2)},\boldsymbol{w}_{0:t-1}^{(2)})\right|\le\beta(\left|x_{0}^{(1)}-x_{0}^{(2)}\right|,t)\nonumber \\
 & +\alpha_{1}(\left\Vert \boldsymbol{w}_{0:t-1}^{(1)}-\boldsymbol{w}_{0:t-1}^{(2)}\right\Vert )+\alpha_{2}(\left\Vert h(\boldsymbol{x}_{0:t}^{(1)})-h(\boldsymbol{x}_{0:t}^{(2)})\right\Vert ),\label{eq:definition - i-IOSS}
\end{align}
where $x_{t}^{(i)}$ is a shorthand of $x_{t}(x_{0}^{(i)},\boldsymbol{w}_{0:t-1}^{(i)})$
for $i=1$ and 2.
\end{defn}
The definition of i-IOSS can be interpreted as a ``detectability''
concept for nonlinear systems \cite{sontag1997output}, as the state
may be ``detected'' from the \emph{noise-free} output by (\ref{eq:definition - i-IOSS}).

In particular, if in \eqref{eq:definition - i-IOSS} $\beta(s,t)=\alpha(s)a^{t}$
for all $s,t\ge0$, with $\alpha\in\K$ and $a$ being a constant
within $(0,1)$, we say that the system is \emph{exponentially i-IOSS}
or \emph{exp-i-IOSS} for short. This can be viewed as extending the
exponential input-to-state stability \cite{grune1999input,liu2010exponential}
to the context of i-IOSS. 
\begin{defn}
\label{def: (-factorizable-} ($\KdL$ function) A $\KL$ function
$\beta$ is called a $\KdL$ function if there exist functions $\alpha\in\K$
and $\varphi\in\LL$ such that $\beta(s,t)=\alpha(s)\varphi(t)$,
for all $s,t\ge0$. 
\end{defn}
As an example, the $\KL$ function $se^{-t}$ is a $\KdL$ function
for $s,t\ge0$. The next lemma shows the general interest of a $\KdL$
function. 
\begin{lem}
\label{lem: factorizable-KL-bound} ($\KdL$ bound, Lemma 8 in \cite{sontag1998comments}) Given an arbitrary
$\KL$ function $\beta$, there exists a $\KdL$ function $\bar{\beta}$
such that $\beta(s,t)\le\bar{\beta}(s,t)$ for all $s,t\ge0$. 
\end{lem}

Lemma \ref{lem: factorizable-KL-bound} implies that the i-IOSS property
defined by means of a $\KL$ function can be defined equivalently
using a $\KdL$ function. This is useful in the later analyses of
FIE and MHE. The following definition of a Lipschitz continuous function
will also be used in the analysis of MHE.

\begin{defn} \label{def: Lipschitz-continuous}
(Lipschitz continuous function) A function $f:\mathbb{R}^{n}\to\mathbb{R}^{m}$
is Lipschitz continuous over a subset $\mathbb{S}\in\mathbb{R}^{n}$
if there is a constant $c$ such that $|f(x)-f(y)|\le c|x-y|$ for
all $x,\thinspace y\in\mathbb{S}$. 
\end{defn}

\section{Optimization-based State Estimation\label{sec: optimization-based estimation}}

Consider a discrete-time nonlinear system described by 
\begin{equation}
x_{t+1}=f(x_{t},w_{t}),\,\,y_{t}=h(x_{t})+v_{t},\label{eq:system}
\end{equation}
where $x_{t}\in\mathbb{R}^{n}$ is the system state, $y_{t}\in\mathbb{R}^{p}$
the measurement, $w_{t}\in\mathbb{R}^{g}$ the process disturbance,
$v_{t}\in\mathbb{R}^{p}$ the measurement disturbance, all at time
$t$. 

{
Since control inputs known up to the estimation time can be treated as given constants, they do not cause difficulty to the later defined optimization and related analyses and hence are ignored for brevity in the problem formulation
\cite{rawlings2012optimization}. 
}
The functions $f$ and $h$ are
assumed to be continuous and known, and the initial state $x_{0}$
and the disturbances $(w_{t},v_{t})$ are modeled as unknown but \emph{bounded}
variables.

Given a time $t$, the state estimation problem is to find
an optimal estimate of state $x_{t}$ based on measurements $\{y_{\tau}\}$
for $\tau$ belonging to a time set and satisfying $\tau\le t$. In the ideal case, all measurements
up to time $t$ are used, leading to the so-called FIE; and in the
practical case, only measurements within a limited distance from time $t$ are used, yielding the so-called MHE. Both FIE
and MHE can be cast as optimization problems as defined next.

Let the decision variables of FIE be $(\boldsymbol{\chi}_{0:t},\boldsymbol{\omega}_{0:t-1},\boldsymbol{\nu}_{0:t})$,
which correspond to the system variables $(\boldsymbol{x}_{0:t},\boldsymbol{w}_{0:t-1},\boldsymbol{v}_{0:t})$,
and the optimal decision variables be $(\hat{\boldsymbol{x}}_{0:t},\hat{\boldsymbol{w}}_{0:t-1},\hat{\boldsymbol{v}}_{0:t})$.
Since $\hat{\boldsymbol{x}}_{0:t}$, which consists of optimal estimates
at all sampled times, is uniquely determined once $\hat{x}_{0}$ and
$\hat{\boldsymbol{w}}_{0:t-1}$ are known, the decision variables
essentially reduce to $(\boldsymbol{\chi}_{0},\boldsymbol{\omega}_{0:t-1},\boldsymbol{\nu}_{0:t})$.

{
Here, although $\boldsymbol{\nu}_{0:t}$ is uniquely determined by $\boldsymbol{\chi}_{0}$
and $\boldsymbol{\omega}_{0:t-1}$, we keep it for the convenience
of expressing bounds and penalty costs to be defined on $\boldsymbol{\nu}_{0:t}$.
}%

Given a present time $t\in\mathbb{I}_{\ge0}$, let $\bar{x}_{0}$
be the prior estimate of the initial state which may be obtained from
the initial or historical measurements. The uncertainty in the initial
state is thus represented by $\chi_{0}-\bar{x}_{0}$. Denote the time-dependent
cost function as $V_{t}(\chi_{0}-\bar{x}_{0},\boldsymbol{\omega}_{0:t-1},\boldsymbol{\nu}_{0:t})$,
which penalizes uncertainties in the initial state, the process and
the measurements. Then, FIE is defined by the following optimization
problem: 
\begin{equation}
\begin{aligned}\mathrm{FIE:}\thinspace\thinspace & \inf V_{t}(\chi_{0}-\bar{x}_{0},\boldsymbol{\omega}_{0:t-1},\boldsymbol{\nu}_{0:t})\\
\mathrm{subject~to,~} & \chi_{\tau+1}=f(\chi_{\tau},\omega_{\tau}),\,\,\forall\tau\in\mathbb{I}_{0:t-1},\\
 & y_{\tau}=h(\chi_{\tau})+\nu_{\tau},\,\,\forall\tau\in\mathbb{I}_{0:t},\\
 & \chi_{0}\in\mathbb{B}_{x}^{0},\thinspace\thinspace\boldsymbol{\omega}_{0:t-1}\in\mathbb{B}_{w}^{0:t-1},\,\,\boldsymbol{\nu}_{0:t}\in\mathbb{B}_{v}^{0:t},
\end{aligned}
\label{eq: FIE}
\end{equation}
where $\{\chi_{0},\thinspace\boldsymbol{\omega}_{0:t-1},\thinspace\boldsymbol{\nu}_{0:t}\}$
are the decision variables. Here $\mathbb{B}_{x}^{0}$, $\mathbb{B}_{w}^{0:t-1}$
and $\mathbb{B}_{v}^{0:t}$ denote the sets of bounded initial states, bounded
sequences of process and measurement disturbances, respectively, for
all $t\in\mathbb{I}_{\ge0}$, of which the latter two sets may vary with time. Since the optimal decision variable $\hat{x}_\tau$, for any $\tau\le t$, is dependent on the time $t$ when the FIE instance is defined, to be unambiguous we use $\hat{x}_{t}^{*}$ to exclusively represent $\hat{x}_{t}$
that is obtained from the FIE instance defined at time $t$. This keeps $\hat{x}_{t}^{*}$
unique, while $\hat{x}_{t}$ varies as the FIE is renewed with new measurements as time elapses.

Given a constant $T\in\mathbb{I}_{\ge0}$, if the measurements are
limited only to the $T+1$ measurements backwards from and including
the present time $t$, then the following optimization defines MHE,
i.e., 
\begin{equation}
\begin{aligned}\mathrm{MHE:~} & \inf V_{T}(\chi_{t-T}-\bar{x}_{t-T},\boldsymbol{\omega}_{t-T:t-1},\boldsymbol{\nu}_{t-T:t})\\
\mathrm{subject~to,~} & \chi_{\tau+1}=f(\chi_{\tau},\omega_{\tau}),\,\,\forall\tau\in\mathbb{I}_{t-T:t-1},\\
 & y_{\tau}=h(\chi_{\tau})+\nu_{\tau},\,\,\forall\tau\in\mathbb{I}_{t-T:t},\\
 \boldsymbol{\chi}_{t-T}\in\mathbb{B}_{x}^{t-T},\,\,& \boldsymbol{\omega}_{t-T:t-1}\in\mathbb{B}_{w}^{t-T:t-1},\,\,\boldsymbol{\nu}_{t-T:t}\in\mathbb{B}_{v}^{t-T:t},
\end{aligned}
\label{eq: MHE}
\end{equation}
where $\{\chi_{t-T},\thinspace\boldsymbol{\omega}_{t-T:t-1},\thinspace\boldsymbol{\nu}_{t-T:t}\}$
are the decision variables, and $\bar{x}_{t-T}$ is a prior estimate
of $x_{t-T}$, and $\mathbb{B}_{x}^{t-T}$, $\mathbb{B}_{w}^{t-T:t-1}$
and $\mathbb{B}_{v}^{t-T:t}$ denote the bounding sets for the time period from $t-T$ to $t$. We use $\hat{x}_{t}^{\star}$ to
represent $\hat{x}_{t}$ that is obtained from the MHE instance defined at time
$t$. By this way, $\hat{x}_{t}^{\star}$ remains unique although
$\hat{x}_{t}$ varies as the MHE instance is renewed with new measurements. Since the cost function is in the same form of FIE
except for the truncated argument variables, the MHE defined in \eqref{eq: MHE}
is named as the \textit{associated} MHE of the FIE defined in \eqref{eq: FIE},
and vice versa. 

{
By definition, FIE uses all available historical measurements to perform state estimation. So its computational complexity increases with time and will ultimately become intractable, which makes FIE impractical for applications. For this reason, FIE is studied mainly for its theoretical interest: its performance can be viewed as
a limit or benchmark that MHE tries to approach, and its stability
can be a good start point for the stability analysis of MHE, which
will become clear later on.
}%

An important issue in designing FIE and MHE is to identify conditions
under which the associated optimizations have optimal solutions such
that the state estimates satisfy the RGAS property defined below.
Let $\boldsymbol{x}_{0:t}(x_{0},\boldsymbol{w}_{0:t-1})$ denote a
state sequence generated from an initial condition $x_{0}$, and a
disturbance sequence $\boldsymbol{w}_{0:t-1}$. 
\begin{defn}
(RGAS \cite{rawlings2012optimization}) The estimate $\hat{x}_{t}$
of the state $x_{t}$ is based on partial or full sequence of the
noisy measurements, $\boldsymbol{y}_{0:t}=h(\boldsymbol{x}_{0:t}(x_{0},\boldsymbol{w}_{0:t-1}))+\boldsymbol{v}_{0:t}$.
The estimate is RGAS if for all $x_{0}$ and $\bar{x}_{0}$, there
exist functions $\beta_{x}\in\KL$ and $\alpha_{w},\alpha_{v}\in\K$
such that the following inequality holds for all $\boldsymbol{\omega}_{0:t-1}\in\mathbb{B}_{w}^{0:t-1}$,
$\boldsymbol{\nu}_{0:t}\in\mathbb{B}_{v}^{0:t}$ and $t\in\mathbb{I}_{\ge0}$:
\begin{align}
 & \left|x_{t}-\hat{x}_{t}\right|\nonumber \\
 & \le\beta_{x}(\left|x_{0}-\bar{x}_{0}\right|,t)+\alpha_{w}(\left\Vert \boldsymbol{w}_{0:t-1}\right\Vert )+\alpha_{v}(\left\Vert \boldsymbol{v}_{0:t}\right\Vert ).\label{def: RGAS}
\end{align}
\end{defn}
Note that, the last measurement $y_{t}$ and hence the corresponding
fitting error $v_{t}$ are considered in the above inequality, which
is however absent in the original definition \cite{rawlings2012optimization}.
To have FIE or MHE that is RGAS, the cost function needs to penalize
the uncertainties appropriately, and meanwhile the system dynamics
should satisfy certain conditions. We present such sufficient conditions
for FIE and MHE respectively in the next two sections.

\begin{rem} \label{rem: on-notation}
In the above formulation of FIE, the notation of an
estimate, $\hat{\cdot}_{\tau}$ for $\tau\in\mathbb{I}_{0:t}$, is
a shorthand of $\hat{\cdot}_{\tau}(x_{0},\boldsymbol{w}_{0:t-1},\boldsymbol{v}_{0:t})$.
Similar shorthand is used in MHE.
The meaning of $\hat{\cdot}_{\tau}$ hence depends on time, and will be explained if ambiguity arises.
\end{rem}

\section{Robust Stability of FIE\label{sec: RGAS of FIE}}

This section summarizes the results on robust stability of FIE
which were obtained in our recent conference paper \cite{hu2015optimization}. The results are rephrased for the ease of
understanding, and some changes are also included and explained. The stability results rely on the following two assumptions.

\begin{assumption} \label{assump: A1} {
The cost function of FIE is given as: $V_{t}(\chi_{0}-\bar{x}_{0},\boldsymbol{\omega}_{0:t-1},\boldsymbol{\nu}_{0:t})=V_{t,1}(\chi_{0}-\bar{x}_{0})+V_{t,2}(\boldsymbol{\omega}_{0:t-1},\boldsymbol{\nu}_{0:t})$, for $t\in\mathbb{I}_{\ge0}$,
where $V_{t,1}$ and $V_{t,2}$ are continuous functions and satisfy
the following inequalities for all $\chi_{0}\in\mathbb{B}_{x}^{0}$, $\boldsymbol{\omega}_{0:t-1}\in\mathbb{B}_{w}^{0:t-1}$,
$\boldsymbol{\nu}_{0:t}\in\mathbb{B}_{v}^{0:t}$ and $t\in\mathbb{I}_{\ge0}$:
\begin{align}
\underbar{\ensuremath{\rho}}_{x}(\left|\chi_{0}-\bar{x}_{0}\right|,t)\le V_{t,1}(\chi_{0}-\bar{x}_{0})\le\rho_{x}(\left|\chi_{0}-\bar{x}_{0}\right|,t),\label{eq: Assumption 1a}\\
\underbar{\ensuremath{\gamma}}_{w}(\left\Vert \boldsymbol{\omega}_{0:t-1}\right\Vert )+\underbar{\ensuremath{\gamma}}_{v}(\left\Vert \boldsymbol{\nu}_{0:t}\right\Vert )\le V_{t,2}(\boldsymbol{\omega}_{0:t-1},\boldsymbol{\nu}_{0:t})\nonumber \\
\le\gamma_{w}(\left\Vert \boldsymbol{\omega}_{0:t-1}\right\Vert )+\gamma_{v}(\left\Vert \boldsymbol{\nu}_{0:t}\right\Vert ),\label{eq: Assumption 1b}
\end{align}
where $\underbar{\ensuremath{\rho}}_{x},\rho_{x}\in\KL$ and $\underbar{\ensuremath{\gamma}}_{w},\underbar{\ensuremath{\gamma}}_{v},\gamma_{w},\gamma_{v}\in\Kinf$.
}%
\end{assumption}

\begin{assumption} \label{assump: A2} The $\K$ and $\KL$ functions
in (\ref{eq:definition - i-IOSS}), (\ref{def: RGAS}) and (\ref{eq: Assumption 1a})
satisfy the following inequalities for all $s_{x},s_{w},s_{v},t\ge0$:
\begin{align}
 & \beta\left(s_{x}+\underbar{\ensuremath{\rho}}_{x}^{-1}\left(\rho_{x}(s_{x},t)+\gamma_{w}(s_{w})+\gamma_{v}(s_{v}),t\right),t\right)\nonumber \\
 & \le\bar{\beta}_{x}(s_{x},t)+\bar{\alpha}_{w}(s_{w})+\bar{\alpha}_{v}(s_{v}), \label{eq: Assumption 2}
\end{align}
in which $\underbar{\ensuremath{\rho}}_{x}^{-1}(\cdot,t)$ refers
to the inverse of $\underbar{\ensuremath{\rho}}_{x}(s,t)$ in its
first argument $s$ at time $t$, and $\bar{\beta}_{x}$, $\bar{\alpha}_{w}$
and $\bar{\alpha}_{v}$ are certain $\KL$, $\K$ and $\K$ functions,
respectively. \end{assumption}

If we interpret the cost function $V_{t}$ as to measure the deviation
of a state estimate of a disturbed system from the true state of
the corresponding undisturbed system, then Assumption \ref{assump: A1} basically
requires that the deviation is both lower and upper bounded by i-IOSS
like limits. 
{
The assumption ensures that the FIE has sufficient sensitivity
to the involved uncertainties, i.e., the sub-cost $V_{t,1}(\chi_{0}-\bar{x}_{0})$ decays neither too fast nor too slowly with respect to $\chi_{0}-\bar{x}_{0}$, and the sub-cost $V_{t,2}(\boldsymbol{\omega}_{0:t-1},\boldsymbol{\nu}_{0:t})$ are lower and upper bounded by strictly increasing functions of $(\boldsymbol{\omega}_{0:t-1},\boldsymbol{\nu}_{0:t})$. 
}%

{
In Assumption
\ref{assump: A2}, if we interpret the argument, $\rho_{x}(s_{x},t)+\gamma_{w}(s_{w})+\gamma_{v}(s_{v})$, as a metric of the deviation from the prior of the true initial state (as caused by the uncertainties $s_x$, $s_w$ and $s_v$), then the $\KL$ function $\underbar{\ensuremath{\rho}}_{x}^{-1}$ aims to infer an upper bound of the deviation based on this metric. The inferred bound, $\underbar{\ensuremath{\rho}}_{x}^{-1}\left(\rho_{x}(s_{x},t)+\gamma_{w}(s_{w})+\gamma_{v}(s_{v}),t\right)$, is added to $s_x$ to form an error bound of the inferred initial state. Subsequently, this error bound is required to be small enough such that the induced error is bounded by an i-IOSS like limit after suppression by the i-IOSS property of the system (as indicated by the $\KL$ function $\beta$ here). Alternatively, we may simply interpret Assumption
\ref{assump: A2} as to require that the FIE is more sensitive than the
system to the uncertainty in the initial state, so that accurate inference of the initial state is possible. The conditions of Assumption \ref{assump: A2} and their interpretations will become more concrete in Lemmas \ref{lem:The-FIE-defined}-\ref{lem:The-FIE-defined3} to be presented ahead.
}

The robust stability of FIE is then established under Assumptions \ref{assump: A1} and \ref{assump: A2}. 
\begin{thm}
\label{thm:RGAS-of-FIE}(RGAS of FIE) The FIE defined in \eqref{eq: FIE}
is RGAS if the three conditions are satisfied: 1) the system described
in \eqref{eq:system} is i-IOSS, 2) the cost function of the FIE satisfies
Assumptions \ref{assump: A1}-\ref{assump: A2}, and 3) the infimum
of the optimization in the FIE is attainable, i.e., exists and is numerically obtainable.
\end{thm}

{
The condition 3) is needed because the state estimate is assumed to be computed as an optimal solution to the optimization defined in \eqref{eq: FIE}.
}%
As a result, the $\K$ and $\KL$ functions of the RGAS
property (cf. Definition \ref{def: RGAS}) can be obtained explicitly
as:
\begin{multline}
\beta_{x}(\left|x_{0}-\bar{x}_{0}\right|,\,t)=\bar{\beta}_{x}\left(\left|x_{0}-\bar{x}_{0}\right|,\,t\right)\\
+\alpha_{1}\left(3\underbar{\ensuremath{\gamma}}_{w}^{-1}\left(3\rho_{x}(\left|x_{0}-\bar{x}_{0}\right|,\,t)\right)\right)\\
+\alpha_{2}\left(3\underbar{\ensuremath{\gamma}}_{v}^{-1}\left(3\rho_{x}(\left|x_{0}-\bar{x}_{0}\right|,\,t)\right)\right),\label{eq: beta_x}
\end{multline}
\begin{multline}
\alpha_{w}(\left\Vert \boldsymbol{w}_{0:t-1}\right\Vert )=\bar{\alpha}_{w}\left(\left\Vert \boldsymbol{w}_{0:t-1}\right\Vert \right)\\
+\alpha_{1}\left(3\left\Vert \boldsymbol{w}_{0:t-1}\right\Vert +3\underbar{\ensuremath{\gamma}}_{w}^{-1}\left(3\gamma{}_{w}(\left\Vert \boldsymbol{w}_{0:t-1}\right\Vert )\right)\right)\\
+\alpha_{2}\left(3\underbar{\ensuremath{\gamma}}_{v}^{-1}\left(3\gamma{}_{w}(\left\Vert \boldsymbol{w}_{0:t-1}\right\Vert )\right)\right),\label{eq: alpha_w}
\end{multline}
\begin{multline}
\alpha_{v}(\left\Vert \boldsymbol{v}_{0:t}\right\Vert )=\bar{\alpha}_{v}\left(\left\Vert \boldsymbol{v}_{0:t}\right\Vert \right)+\alpha_{1}\left(3\underbar{\ensuremath{\gamma}}_{w}^{-1}\left(3\gamma{}_{v}(\left\Vert \boldsymbol{v}_{0:t}\right\Vert )\right)\right)\\
+\alpha_{2}\left(3\left\Vert \boldsymbol{v}_{0:t}\right\Vert +3\underbar{\ensuremath{\gamma}}_{v}^{-1}\left(3\gamma{}_{v}(\left\Vert \boldsymbol{v}_{0:t}\right\Vert )\right)\right).\label{eq: alpha_v}
\end{multline}

\begin{rem}
As shown in \cite{hu2015optimization}, FIE can prove to converge
to the true state when the disturbances are convergent to zero if
the feasible sets $\mathbb{B}_{w}^{0:t-1}$ and $\mathbb{B}_{v}^{0:t}$ restrict
the disturbances estimates to be convergent to zero. However, it is
unclear if there exists a form of cost function such that the conclusion
remains true without imposing this restriction. 
\end{rem}

A more specific form of the sub-cost function $V_{t,2}$ that satisfies
Assumption \ref{assump: A1} is given by the following:{\small{}{} 
\begin{align}
 & V_{t,2}(\boldsymbol{\omega}_{0:t-1},\,\boldsymbol{\nu}_{0:t}):=\dfrac{\lambda_{w}}{t}\sum_{\tau\in\mathbb{I}_{0:t-1}}l_{w,\tau}(\omega_{\tau})+\dfrac{\lambda_{v}}{t+1}\sum_{\tau\in\mathbb{I}_{0:t}}l_{v,\tau}(\nu_{\tau})\label{eq: specific form of the sub-cost}\\
 & \,\,\,\,+(1-\lambda_{w})\max_{\tau\in\mathbb{I}_{0:t-1}}l_{w,\tau}(\omega_{\tau})+(1-\lambda_{v})\max_{\tau\in\mathbb{I}_{0:t}}l_{w,\tau}(\nu_{\tau}),\nonumber 
\end{align}
}for given constants $\lambda_{w},\,\lambda_{v}\in[0,\,1)$, in which
the functions $l_{w,\tau}$ and $l_{v,\tau}$ satisfy the following
inequalities for all $\boldsymbol{\omega}_{0:t-1}\in\mathbb{B}_{w}^{0:t-1}$
and $\boldsymbol{\nu}_{0:t}\in\mathbb{B}_{v}^{0:t}$: 
\begin{equation}
\begin{gathered}\underbar{\ensuremath{\gamma}}_{w}'(|\omega_{\tau}|)\le l_{w,\tau}(\omega_{\tau})\le\gamma_{w}'(|\omega_{\tau}|),\\
\underbar{\ensuremath{\gamma}}_{v}'(|\nu_{\tau}|)\le l_{v,\tau}(\nu_{\tau})\le\gamma_{v}'(|\nu_{\tau}|),
\end{gathered}
\label{eq: specific-subcost-form - b}
\end{equation}
where $\underbar{\ensuremath{\gamma}}_{w}',\,\underbar{\ensuremath{\gamma}}_{v}',\,\gamma_{w}',\,\gamma_{v}'\in\Kinf$.
In the sub-cost function, the terms associated with $\omega_{\tau}$
vanish if $t=0$.

On the other hand, more specific forms of the sub-cost function $V_{t,1}$ that satisfies Assumptions \ref{assump: A1} and \ref{assump: A2} can be obtained if the $\KL$ function $\beta$ of the i-IOSS property
of the system belongs to one of two particular types. The derivation is based
on the next lemma.
\begin{lem}
\label{lem:The-FIE-defined}Assumption \ref{assump: A2} is satisfied
if the $\KL$ functions $\beta$ in \eqref{eq:definition - i-IOSS}
and $\underbar{\ensuremath{\rho}}_{x},\thinspace\rho_{x}$ in \eqref{eq: Assumption 1a}
are $\KdL$ functions in the form of $\beta(s,\,t)=\mu_{1}(s)\varphi_{1}(t)$,
$\underbar{\ensuremath{\rho}}_{x}(s,\,t)=\mu_{2}(s)\varphi_{2}(t)$
and $\ensuremath{\rho}_{x}(s,\,t)=\mu_{3}(s)\varphi_{2}(t)$, with
$\mu_{1},\,\mu_{2},\thinspace\mu_{3}\in\K$ and $\varphi_{1},\varphi_{2}\in\LL$,
and further for any $\pi\in\K$ there exists $\pi'\in\K$ such that
\begin{equation}
\mu_{1}\left(3\mu_{2}^{-1}\left(\frac{\pi(s)}{\varphi_{2}(t)}\right)\right)\varphi_{1}(t)\le\pi'(s)\label{eq: key-stability-condition}
\end{equation}
which holds for all $s,\thinspace t\ge0$. 
\end{lem}
\begin{IEEEproof}
The proof is the same as that of Corollary 1 in \cite{hu2015optimization}
except that a tighter upper bound is used during the deduction:
\begin{align}
 & \beta\left(s_{x}+\underbar{\ensuremath{\rho}}_{x}^{-1}\left(\rho_{x}(s_{x},\,t)+\gamma_{w}(s_{w})+\gamma_{v}(s_{v}),\, t\right),\,t\right)\nonumber \\
 & \le\beta\left(3s_{x}+3\underbar{\ensuremath{\rho}}_{x}^{-1}\left(3\rho_{x}(s_{x},\,t),\,t\right),\,t\right)\nonumber \\
 & \,\,\,\,\,\,+\beta\left(3\underbar{\ensuremath{\rho}}_{x}^{-1}\left(3\gamma_{w}(s_{w}),\,t\right),\,t\right)+\beta\left(3\underbar{\ensuremath{\rho}}_{x}^{-1}\left(3\gamma_{v}(s_{v}),\,t\right),\,t\right),\nonumber
\end{align}
of which the three terms can be proved to be upper bounded by $\KL$, $\K$ and $\K$ functions, respectively, by use of \eqref{eq: key-stability-condition}. 
\end{IEEEproof}
In Lemma \ref{lem:The-FIE-defined}, the assumption of $\beta$ being
a $\KdL$ function is trivial because it is always feasible to assign
such a function as an alternative if the original $\KL$ function
$\beta$ is not in a $\KdL$ form (cf. Lemma \ref{lem: factorizable-KL-bound}).
{
The condition that $\underbar{\ensuremath{\rho}}_{x}$ and $\rho_{x}$
in \eqref{eq: Assumption 1a} are $\KdL$ functions is not imposed on the
system dynamics, but a requirement on the cost function of
FIE. 
}%
The key condition thus boils down to \eqref{eq: key-stability-condition},
which is basically an alternative of the more general condition \eqref{eq: Assumption 2}. Therefore the previous interpretation of \eqref{eq: Assumption 2} (or Assumption \ref{assump: A2}) is applicable to \eqref{eq: key-stability-condition}.

Based on Lemma \ref{lem:The-FIE-defined}, we can prove that the FIE
admits a even more specific cost function if the system is i-IOSS
with a $\KL$ bound in the rational form. 
\begin{lem}
\label{lem:The-FIE-defined2}Assumption \ref{assump: A2} is satisfied
if the three conditions are satisfied: a) the system \eqref{eq:system}
is i-IOSS as per \eqref{eq:definition - i-IOSS} in which the $\KL$
bound is explicitly given as $\beta(s,t)=c_{1}s^{a_{1}}(t+1)^{-b_{1}}$
for some constants $c_{1},a_{1},b_{1}>0$ and all $s,t\ge0$, and
b) the sub-cost function $V_{t,1}$ is defined as 
\begin{align*}
V_{t,1}(\mathcal{X}_{0}-\bar{x}_{0}) & =c_{2}|\mathcal{X}_{0}-\bar{x}_{0}|^{a_{2}}(t+1)^{-b_{2}},
\end{align*}
with $a_{2},\thinspace b_{2}>0$, and c) the parameters $a_{2}$ and
$b_{2}$ satisfy $\frac{a_{2}}{b_{2}}\ge\frac{a_{1}}{b_{1}}$. 
\end{lem}
This lemma implies the main result of \cite{ji2016robust} if the
design parameter $b_{2}$ is fixed to 1 (with a minor difference that
here the FIE is able to utilize the last measurement in the estimation,
whose fitting error is penalized through $\nu(t)$).

Moreover, if the system described in \eqref{eq:system} is exp-i-IOSS,
then the conclusion remains valid by replacing the rational form of $\KL$
bound in Lemma \ref{lem:The-FIE-defined2} with an exponential form. 
\begin{lem}
\label{lem:The-FIE-defined3}Assumption \ref{assump: A2} is satisfied
if the three conditions are satisfied: a) the system \eqref{eq:system}
is exp-i-IOSS as per \eqref{eq:definition - i-IOSS} in which the
$\KL$ function is explicitly given as $\beta(s,\,t)=c_{1}s^{a_{1}}b_{1}^{t}$
for some constants $c_{1},a_{1}>0$ and $0<b_{1}<1$ and all $s,t\ge0$,
and b) the sub-cost function $V_{t,1}$ is defined as 
\begin{align*}
V_{t,1}(\mathcal{X}_{0}-\bar{x}_{0}) & =c_{2}|\mathcal{X}_{0}-\bar{x}_{0}|^{a_{2}}b_{2}^{t},
\end{align*}
with $a_{2}>0$ and $0<b_{2}<1$, and c) the parameters $a_{2}$ and
$b_{2}$ satisfy $\sqrt[a_{2}]{b}_{2}\ge\sqrt[a_{1}]{b}_{1}$. 
\end{lem}
In condition c) of Lemma \ref{lem:The-FIE-defined3}, the constraint
$b_{2}<1$ is required to make sure that $c_{2}s^{a_{2}}b_{2}^{t}$
is a $\KL$ function of $s$ and $t$, so as to satisfy Assumption
\ref{assump: A1} of Theorem \ref{thm:RGAS-of-FIE}.
\begin{rem}
As shown in \cite{glas1987exponential}, the set of exponentially
stable systems are dense in the whole set of asymptotically stable
systems. So it seems not to lose generality to assume exp-i-IOSS
systems in practice as in Lemma \ref{lem:The-FIE-defined3}. 
\end{rem}

\section{Robust Stability of MHE\label{sec: RGAS of the MHE}}

{
At any discrete time, an MHE instance can be treated as an associated FIE instance that is confined to the same optimization horizon. Thus, the associated FIE instance being RGAS
implies the RGAS of MHE within its present optimization horizon.
If we interpret this as MHE being robust \textit{locally} asymptotically stable (RLAS) within each optimization horizon of a given size,
then the challenge reduces to identifying the conditions under which
RLAS implies RGAS of MHE.
}

To that end, we need an assumption on the prior estimate of the initial
state of each MHE instance.

\begin{assumption} \label{assump: A3} Given any time $t\ge T+1$,
the prior estimate $\bar{x}_{t-T}$ of $x_{t-T}$ satisfies the following
constraint: 
\[
|x_{t-T}-\bar{x}_{t-T}|\le|x_{t-T}-\hat{x}_{t-T}^{\star}|.
\]
\end{assumption}

The assumption is trivially satisfied if $\bar{x}_{t-T}$ is set to
$\hat{x}_{t-T}^{\star}$, which is the past MHE estimate obtained at time
$t-T$. Alternatively, a better $\bar{x}_{t-T}$ might be obtained
with smoothing techniques which use measurements both before and after
time $t-T$ \cite{aravkin2016generalized}. Since a rigorous derivation is non-trivial, the extension is left for future research.

The next lemma links the
robust stability and convergence of MHE with those of its associated
FIE. 
\begin{lem}\label{lem: MHE-FIE}
(Stability link from FIE to MHE) Consider the
MHE under Assumption \ref{assump: A3}. Let the uncertainty in the
initial state be bounded as $|x_{0}-\bar{x}_{0}|\le M_{0}$, and the
disturbances be bounded as $|w_{t}|\le M_{w}$ and $|v_{t}|\le M_{v}$
for all $t\in\mathbb{I}_{\ge0}$. Given a constant $\eta\in(0,\thinspace1)$,
the following two conclusions hold:

{
a) If the associated FIE is RGAS as per \eqref{def: RGAS}, in which
the $\KL$ function $\beta_{x}$ satisfies $\beta_{x}(s,\thinspace t)\le\mu(s)\varphi(t)$
for some $\mu\in\K$, $\varphi\in\LL$ and all $s\in\mathbb{R}_{\ge0}$,
$t\in\mathbb{I}_{\ge0}$ , and if there exists $T_{\bar{s},\eta}$
such that $\mu(s)\varphi(T_{\bar{s},\eta})\le\eta s$ for all $s\in[0,\thinspace\bar{s}]$
with $\bar{s}:=\beta_{x}(M_{0},\thinspace0)+\frac{1}{1-\eta}\left(\alpha_{w}(M_{w})+\alpha_{v}(M_{v})\right)$,
then MHE is RGAS for all $T\ge T_{\bar{s},\eta}$. In particular,
if $\mu(s)$ is Lipschitz continuous at the origin, then $T_{\bar{s},\eta}$
exists and can be determined from the inequality, $\varphi(T_{\bar{s},\eta})\le\frac{\eta s^{*}}{\mu(s^{*})}$
with $s^{*}:=\arg\min_{s\in[0,\thinspace\bar{s}]}\frac{s}{\mu(s)}$.
}%

b) If the associated FIE estimate ($\hat{x}_{t}^{*}$) converges to
the true state, i.e., the estimate satisfies $|x_{t}-\hat{x}_{t}^{*}|\le\rho'(|x_{0}-\bar{x}_{0}|,\thinspace t)\le\mu'(|x_{0}-\bar{x}_{0}|)\varphi'(t)$
for some $\rho'\in\KL$, $\mu'\in\K$, $\varphi'\in\LL$ and all $t\in\mathbb{I}_{\ge0}$,
and if there exists $T_{\bar{s}',\eta}$ such that $\mu(s)\varphi(T_{\bar{s}',\eta})\le\eta s$
for all $s\in[0,\thinspace\bar{s}']$ with $\bar{s}':=\rho'(M_{0},\thinspace0)$,
then the MHE estimate ($\hat{x}_{t}^{\star}$) converges to the true
state for all $T\ge T_{\bar{s}',\eta}$. In particular, if $\mu'(s)$
is Lipschitz continuous at the origin, then $T_{\bar{s}',\eta}$ exists
and can be determined from the inequality, $\varphi(T_{\bar{s}',\eta})\le\frac{\eta s^{\star}}{\mu(s^{\star})}$
with $s^{\star}:=\arg\min_{s\in[0,\thinspace\bar{s}']}\frac{s}{\mu'(s)}$.
\end{lem}
\begin{IEEEproof}
a)\emph{ RGAS.} Given $t\in\mathbb{I}_{0:T-1}$, the MHE estimate
($\hat{x}_{t}^{\star}$) is the same as the associated FIE estimate
($\hat{x}_{t}^{*}$), and so the estimation error norm $\left|x_{t}-\hat{x}_{t}^{\star}\right|$
satisfies the RGAS inequality by \eqref{def: RGAS}. Specifically,
under Assumption \ref{assump: A3} the RGAS inequality implies that,
for all $t\in\mathbb{I}_{0:T-1}$,
\begin{equation}
\left|x_{t}-\bar{x}_{t}\right|\le\beta_{x}(M_{0},\thinspace0)+\alpha_{w}(M_{w})+\alpha_{v}(M_{v})\le\bar{s}.\label{eq: error bound from 0 to T}
\end{equation}
Next, we proceed to prove that the RGAS property is maintained for
all $t\in\mathbb{I}_{\ge T}$. 

Given $t\in\mathbb{I}_{\ge T}$, define $n=\left\lfloor \frac{t}{T}\right\rfloor $,
which is the largest integer that is less than or equal to $\frac{t}{T}$.
So $t-nT$ belongs to the set $\mathbb{I}_{0:T-1}$, and hence $\left|x_{t-nT}-\bar{x}_{t-nT}\right|$
satisfies the preceding inequality, i.e., $\left|x_{t-nT}-\bar{x}_{t-nT}\right|\le\bar{s}$.
Treat the MHE defined at time $t-(n-1)T$ as the associated FIE confined
to the time interval $[t-nT,\thinspace t-(n-1)T]$. Therefore, the MHE
satisfies the RGAS property within this interval, that is, by \eqref{def: RGAS}
we have:
\begin{align*}
 & \left|x_{t-(n-1)T}-\hat{x}_{t-(n-1)T}^{\star}\right|\le\beta_{x}(\left|x_{t-nT}-\bar{x}_{t-nT}\right|,\thinspace T)\\
 & +\alpha_{w}(\left\Vert \boldsymbol{w}_{t-nT:t-(n-1)T-1}\right\Vert )+\alpha_{v}(\left\Vert \boldsymbol{v}_{t-nT:t-(n-1)T}\right\Vert )\\
 & \le\mu(\left|x_{t-nT}-\bar{x}_{t-nT}\right|)\varphi(T)+\alpha_{w}(\left\Vert \boldsymbol{w}_{0:t-1}\right\Vert )+\alpha_{v}(\left\Vert \boldsymbol{v}_{0:t}\right\Vert ).
\end{align*}
Since $\left|x_{t-nT}-\bar{x}_{t-nT}\right|\in[0,\thinspace\bar{s}]$
and $\varphi(T)$ decreases with $T$, for all $T\ge T_{\bar{s},\eta}$
we have
\begin{align*}
 & \left|x_{t-(n-1)T}-\hat{x}_{t-(n-1)T}^{\star}\right|\\
 & \le\mu(\left|x_{t-nT}-\bar{x}_{t-nT}\right|)\varphi(T_{\bar{s},\eta})+\alpha_{w}(\left\Vert \boldsymbol{w}_{0:t-1}\right\Vert )+\alpha_{v}(\left\Vert \boldsymbol{v}_{0:t}\right\Vert )\\
 & \le\eta\left|x_{t-nT}-\bar{x}_{t-nT}\right|+\alpha_{w}(\left\Vert \boldsymbol{w}_{0:t-1}\right\Vert )+\alpha_{v}(\left\Vert \boldsymbol{v}_{0:t}\right\Vert )
\end{align*}
where the second inequality follows from the definition of $T_{\bar{s},\eta}$.
Repeat the above reasoning for the MHE defined at time $t-(n-2)T$
with $T\ge T_{\bar{s},\eta}$, yielding
\begin{align*}
 & \left|x_{t-(n-2)T}-\hat{x}_{t-(n-2)T}^{\star}\right|\\
 & \le\beta_{x}(\left|x_{t-(n-1)T}-\bar{x}_{t-(n-1)T}\right|,\thinspace T)\\
 & \quad+\alpha_{w}(\left\Vert \boldsymbol{w}_{t-(n-1)T:t-(n-2)T-1}\right\Vert )\\
 & \quad+\alpha_{v}(\left\Vert \boldsymbol{v}_{t-(n-1)T:t-(n-2)T}\right\Vert )\\
 & \overset{\text{Assump. \ref{assump: A3}}}{\le}\beta_{x}(\left|x_{t-(n-1)T}-\hat{x}_{t-(n-1)T}^{\star}\right|,\thinspace T)\\
 & \quad+\alpha_{w}(\left\Vert \boldsymbol{w}_{0:t-1}\right\Vert )+\alpha_{v}(\left\Vert \boldsymbol{v}_{0:t}\right\Vert )\\
 & \le\beta_{x}\left(\eta\left|x_{t-nT}-\bar{x}_{t-nT}\right|+\alpha_{w}(\left\Vert \boldsymbol{w}_{0:t-1}\right\Vert )+\alpha_{v}(\left\Vert \boldsymbol{v}_{0:t}\right\Vert ),\thinspace T\right)\\
 & \quad+\alpha_{w}(\left\Vert \boldsymbol{w}_{0:t-1}\right\Vert )+\alpha_{v}(\left\Vert \boldsymbol{v}_{0:t}\right\Vert )\\
 & \le\mu\left(\eta\left|x_{t-nT}-\bar{x}_{t-nT}\right|+\alpha_{w}(\left\Vert \boldsymbol{w}_{0:t-1}\right\Vert )+\alpha_{v}(\left\Vert \boldsymbol{v}_{0:t}\right\Vert )\right)\\
 & \quad\times\varphi(T_{\bar{s},\eta})+\alpha_{w}(\left\Vert \boldsymbol{w}_{0:t-1}\right\Vert )+\alpha_{v}(\left\Vert \boldsymbol{v}_{0:t}\right\Vert )\\
 & \le\eta\left(\eta\left|x_{t-nT}-\bar{x}_{t-nT}\right|+\alpha_{w}(\left\Vert \boldsymbol{w}_{0:t-1}\right\Vert )+\alpha_{v}(\left\Vert \boldsymbol{v}_{0:t}\right\Vert )\right)\\
 & \quad+\alpha_{w}(\left\Vert \boldsymbol{w}_{0:t-1}\right\Vert )+\alpha_{v}(\left\Vert \boldsymbol{v}_{0:t}\right\Vert )\\
 & =\eta^{2}\left|x_{t-nT}-\bar{x}_{t-nT}\right| \\
 & \quad+(1+\eta)\left(\alpha_{w}(\left\Vert \boldsymbol{w}_{0:t-1}\right\Vert )+\alpha_{v}(\left\Vert \boldsymbol{v}_{0:t}\right\Vert )\right).
\end{align*}
In the deduction, the inequality \eqref{eq: error bound from 0 to T}
has been used to show that $\eta\left|x_{t-nT}-\bar{x}_{t-nT}\right|+\alpha_{w}(\left\Vert \boldsymbol{w}_{0:t-1}\right\Vert )+\alpha_{v}(\left\Vert \boldsymbol{v}_{0:t}\right\Vert )\le\bar{s}$, and
so the inequality $\mu(s)\varphi(T_{\bar{s},\eta})\le\eta s$
remains applicable.

By induction, we obtain
\begin{align*}
 \left|x_{t}-\hat{x}_{t}^{\star}\right| & \le\eta^{n}\left|x_{t-nT}-\bar{x}_{t-nT}\right| \\
 & \quad+\sum_{i=0}^{n-1}\eta^{i}\left(\alpha_{w}(\left\Vert \boldsymbol{w}_{0:t-1}\right\Vert )+\alpha_{v}(\left\Vert \boldsymbol{v}_{0:t}\right\Vert )\right)\\
 & \le\eta^{\left\lfloor \frac{t}{T}\right\rfloor }\left|x_{t-nT}-\bar{x}_{t-nT}\right| \\
 & \quad+\frac{1}{1-\eta}\left(\alpha_{w}(\left\Vert \boldsymbol{w}_{0:t-1}\right\Vert )+\alpha_{v}(\left\Vert \boldsymbol{v}_{0:t}\right\Vert )\right),
\end{align*}
for all $T\ge T_{\bar{s},\eta}$. Since MHE satisfies the RGAS property
within the time interval $[0,\thinspace t-nT]$, it follows that
\begin{align*}
 & \left|x_{t-nT}-\bar{x}_{t-nT}\right|\le x_{t-nT}-\hat{x}_{t-nT}^{\star}\\
 & \le\beta_{x}(\left|x_{0}-\bar{x}_{0}\right|,\thinspace t-nT)\\
 & \quad+\alpha_{w}(\left\Vert \boldsymbol{w}_{0:t-nT-1}\right\Vert )+\alpha_{v}(\left\Vert \boldsymbol{v}_{0:t-nT}\right\Vert )\\
 & \le\beta_{x}(\left|x_{0}-\bar{x}_{0}\right|,\thinspace0)+\alpha_{w}(\left\Vert \boldsymbol{w}_{0:t-1}\right\Vert )+\alpha_{v}(\left\Vert \boldsymbol{v}_{0:t}\right\Vert ).
\end{align*}
Consequently, 
\begin{align*}
& \left|x_{t}-\hat{x}_{t}^{\star}\right| \\
& \le\eta^{\left\lfloor \frac{t}{T}\right\rfloor }\left(\beta_{x}(\left|x_{0}-\bar{x}_{0}\right|,\thinspace0)+\alpha_{w}(\left\Vert \boldsymbol{w}_{0:t-1}\right\Vert )+\alpha_{v}(\left\Vert \boldsymbol{v}_{0:t}\right\Vert )\right)\\
 & \quad+\frac{1}{1-\eta}\left(\alpha_{w}(\left\Vert \boldsymbol{w}_{0:t-1}\right\Vert )+\alpha_{v}(\left\Vert \boldsymbol{v}_{0:t}\right\Vert )\right)\\
 & \le\eta^{\left\lfloor \frac{t}{T}\right\rfloor }\beta_{x}(\left|x_{0}-\bar{x}_{0}\right|,\thinspace0) \\
 & \quad+\frac{2-\eta}{1-\eta}\left(\alpha_{w}(\left\Vert \boldsymbol{w}_{0:t-1}\right\Vert )+\alpha_{v}(\left\Vert \boldsymbol{v}_{0:t}\right\Vert )\right)\\
 & =:\beta_{x}'(\left|x_{0}-\bar{x}_{0}\right|,\thinspace t)+\alpha_{w}'(\left\Vert \boldsymbol{w}_{0:t-1}\right\Vert )+\alpha_{v}'(\left\Vert \boldsymbol{v}_{0:t}\right\Vert ),
\end{align*}
for all $T\ge T_{\bar{s},\eta}$ , where $\beta_{x}'\in\KL$ and $\alpha_{w}',\thinspace\alpha_{v}'\in\K$.
Therefore, the MHE satisfies the RGAS property for all $t\in\mathbb{I}_{\ge0}$,
which completes the proof of the major conclusion.

If $\mu(s)$ is Lipschitz continuous at the origin, together with
the property that $\mu(0)=0$ and $\mu(s)$ is non-negative and strictly
increasing for all $s\in\mathbb{R}_{\ge0}$, it follows that the value
of $\frac{\mu(s)}{s}$ must be positive and bounded above for all
$s\in[0,\thinspace\bar{s}]$. Consequently, $\frac{s}{\mu(s)}$ is
positive and bounded below. That is, the minimizer $s^{*}:=\arg\min_{s\in[0,\thinspace\bar{s}]}\frac{s}{\mu(s)}$
exists and is well-defined. By the property of a $\LL$ function,
it follows that there exists $T_{\bar{s},\eta}>0$ such that $\varphi(T_{\bar{s},\eta})\le\frac{\eta s^{*}}{\mu(s^{*})}$
and the MHE is RGAS. This proves the rest part of the conclusion.

b)\emph{ Convergence.} If the associated FIE estimate ($\hat{x}_{t}^{*}$)
converges to the true state, then by Lemma 4.5 of \cite{khalil2002nonlinear}
and Lemma \ref{lem: factorizable-KL-bound} in Section \ref{sec: Notation-and-Preliminaries},
there exist $\rho'\in\KL$, $\mu'\in\K$ and $\varphi'\in\LL$ such
that 
\[
|x_{t}-\hat{x}_{t}^{*}|\le\rho'(|x_{0}-\bar{x}_{0}|,\thinspace t)\le\mu'(|x_{0}-\bar{x}_{0}|)\varphi'(t),
\]
for all $t\in\mathbb{I}_{\ge0}$. Continue the proof per part a) but
with the $\KL$ function $\beta_{x}(s,\thinspace t)$ replaced with
$\rho'(s,\thinspace t)$ and the $\K$ functions $\alpha_{w}$ and
$\alpha_{v}$ set to zero. We reach the conclusion that $|x_{t}-\hat{x}_{t}^{\star}|\le\varrho(|x_{0}-\bar{x}_{0}|,\thinspace t)$
for some $\varrho\in\KL$, if $T\ge T_{\bar{s}',\eta}$. This implies
that the MHE estimate ($\hat{x}_{t}^{\star}$) converges to the true
state ($x_{t}$), which completes the proof of the major conclusion.
The rest of the proof with $\mu'(s)$ being Lipschitz continuous at
the origin is completed as per the last paragraph of the proof in
part a).
\end{IEEEproof}

{
From the above proof, we see that $\bar{s}$ and $\bar{s}'$ in Lemma \ref{lem: MHE-FIE} are basically the upper bounds of the uncertainty in the initial state of an MHE instance defined at any time. They are used to define the ranges of the uncertainty within which the conditions of the lemma need to hold. This avoids a stronger condition which assumes $\bar{s}$ or $\bar{s}'$ to be infinite.
}

Lemma \ref{lem: MHE-FIE} indicates that the robust stability of MHE
is implied by the \textit{enhanced} robust stability of its associated
FIE. The enhancing condition requires the moving horizon size ($T$)
to be large enough such that the inequality, $\mu(s)\varphi(T)\le\eta s$,
holds true when the initial state estimation error ($s$) takes a
value within a bounded range. (The lower bound on $T$ can be less
conservative if the size of the moving horizon adapts to the variable
$s$ while keeping the inequality satisfied.) 
{
With $0<\eta<1$, the
condition basically requires each MHE instance to be based on sufficient measurements
so that the effect of the estimation error of the initial state decays over time.
}

The conditions of Lemma \ref{lem: MHE-FIE} become more specific if
the $\K$ function $\mu$ and $\mu'$ have special forms. For example,
if $\mu(s)=\mu'(s)=c_{1}s^{a_{1}}$ with $a_{1}\ge1$ and $c_{1}>0$,
both of which are Lipschitz continuous at the origin, then the conditions
of conclusion a) reduce to that $T\ge T_{\bar{s},\eta}$, satisfying
$\varphi(T_{\bar{s},\eta})\le\frac{\eta}{c_{1}\bar{s}^{a_{1}-1}}$,
which further degenerates to $\varphi(T_{\bar{s},\eta})\le\frac{\eta}{c_{1}}$
if $a_{1}=1$, and meanwhile the conditions of conclusion b) reduce
to $T\ge T_{\bar{s}',\eta}$, satisfying $\varphi'(T_{\bar{s}',\eta})\le\frac{\eta}{c_{1}\bar{s}'^{a_{1}-1}}$,
which further degenerates to $\varphi'(T_{\bar{s}',\eta})\le\frac{\eta}{c_{1}}$
if $a_{1}=1$. Note that, if $0<a_{1}<1$ in these two cases, then
$\mu(s)$ and $\mu'(s)$ are not Lipschitz continuous at the origin
and the RGAS of MHE may not follow.

With the explicit link established between the stability of MHE and
that of its associated FIE, we are able to prove the RGAS of MHE by
enhancing the conditions that establish the RGAS of FIE. In the following,
the symbol $\bar{s}$ remains to be the constant defined in Lemma
\ref{lem: MHE-FIE}. 
\begin{thm}
\label{thm: RGAS-of-MHE}(RGAS of MHE) Suppose that the system described
in \eqref{eq:system} is i-IOSS and the infimum of the MHE defined
in \eqref{eq: MHE} is attainable (i.e., exists and numerically obtainable). Given Assumption \ref{assump: A3}
and any $\eta\in(0,\thinspace1)$, the MHE is RGAS for all $T\ge T_{\eta,\bar{s}}$
if its associated FIE satisfies Assumptions \ref{assump: A1}-\ref{assump: A2}
and the involved $\K$ and $\KL$ functions satisfy 
\begin{multline}
\bar{\beta}_{x}(s,\thinspace T_{\eta,\bar{s}})+\alpha_{1}\left(3\underbar{\ensuremath{\gamma}}_{w}^{-1}\left(3\rho_{x}(s,\thinspace T_{\eta,\bar{s}})\right)\right)\\
+\alpha_{2}\left(3\underbar{\ensuremath{\gamma}}_{v}^{-1}\left(3\rho_{x}(s,\thinspace T_{\eta,\bar{s}})\right)\right)\le\eta s,\label{eq: RGAS-of-MHE-extrac-condition}
\end{multline}
for all $s\in[0,\thinspace\bar{s}]$ and $t\in\mathbb{I}_{\ge0}$. Furthermore, if both disturbance $w(t)$ and noise $v(t)$ converge to zero as $t$ goes to infinity, then the MHE estimate $x^{\star}(t)$ converges to the true state $x(t)$.
\end{thm}
\begin{IEEEproof}
(a) \textit{RGAS}. Under the conditions excluding Assumption \ref{assump: A3} and inequality
\eqref{eq: RGAS-of-MHE-extrac-condition}, the FIE associated with
the MHE is RGAS by Theorem \ref{thm:RGAS-of-FIE}. In the resulting RGAS property,
the $\KL$ bound function is obtained as $\beta_{x}\left(s,\,t\right)=\bar{\beta}_{x}\left(s,\,t\right)+\alpha_{1}\left(3\underbar{\ensuremath{\gamma}}_{w}^{-1}\left(3\rho_{x}(s,\,t)\right)\right)+\alpha_{2}\left(3\underbar{\ensuremath{\gamma}}_{v}^{-1}\left(3\rho_{x}(s,\,t)\right)\right)$
(cf. \eqref{eq: beta_x}). Then, inequality \eqref{eq: RGAS-of-MHE-extrac-condition}
implies that $\beta_{x}\left(s,\,T\right)\le\eta s$ for all $T\ge T_{\eta,\bar{s}}$
and $s\in[0,\thinspace\bar{s}]$. Consequently the conclusion follows
from conclusion a) of Lemma \ref{lem: MHE-FIE}. 

(b) \textit{Convergence}. Since the disturbance $w_t$ and the noise
$v_t$ converge to zero, for any $\epsilon>0$, there exists a time
$t_{\epsilon}$ such that $|w_t|<\epsilon/3$ and $|v_t|<\epsilon/3$
for all $t\ge t_{\epsilon}$. By the definition of $\KL$ function,
there also exists a time $\tau_{\epsilon}$ such that $\beta_{x}(|x_{t_{\epsilon}}-\bar{x}_{t_{\epsilon}}|,\thinspace\tau_{\epsilon})<\epsilon/3$.
Given $t\ge t_{\epsilon}+\tau_{\epsilon}$, by the induction in part (a) of the proof
of Lemma \ref{lem: MHE-FIE}, we observe that, under the conditions
of this theorem, the same RGAS inequality \eqref{def: RGAS} of the MHE
remains valid if the intermediate state $x_{t_{\epsilon}}$ is treated
as the initial state. Consequently, we obtain
\begin{align*}
|\hat{x}^{\star}_{t}-x_{t}| & \le \beta_{x}(|x_{t_{\epsilon}}-\bar{x}_{t_{\epsilon}}|,\thinspace t-t_{\epsilon})\\
 & \quad+\alpha_{w}(\left\Vert \boldsymbol{w}_{t_{\epsilon}:t-1}\right\Vert )+\alpha_{v}(\left\Vert \boldsymbol{v}_{t_{\epsilon}:t}\right\Vert )\\
 & <\epsilon/3+\epsilon/3+\epsilon/3=\epsilon,
\end{align*}
which implies that the MHE estimate $\hat{x}^{\star}_{t}$ converges to $x_{t}$
as $t$ goes to infinity. This completes the proof.
\end{IEEEproof}

{
The proof shows that the left-hand-side of \eqref{eq: RGAS-of-MHE-extrac-condition} is nothing but the $\KL$ bound component in the RGAS property of the FIE that is associated with the MHE. Inequality \eqref{eq: RGAS-of-MHE-extrac-condition} basically requires the $\KL$ bound to be contractive with respect to the estimation error of the initial state for each MHE instance. This is made possible by requiring the MHE to implement a sufficiently large moving horizon as indicated by $T_{\eta,\bar{s}}$.
}

\begin{rem} \label{rem: convergence-of-FIE}
The convergence of FIE is not proved for the same
conditions given in Theorem \ref{thm: RGAS-of-MHE}. This is because, given an initial condition ($x_{0}$), 
the RGAS property of FIE is exclusively associated with $x_0$ and is not applicable if an intermediate state ($x_{t_{\epsilon}}$, $\forall t_{\epsilon}\in\mathbb{I}_{>0}$) is used to replace $x_0$ in the property. This makes it invalid to apply a similar RGAS inequality to establish the convergence as per the above proof.
\end{rem}

\begin{lem}
\label{lem: RGAS of MHE - KdL form}Given conditions a)-c) of Lemma
\ref{lem:The-FIE-defined}, the condition \eqref{eq: RGAS-of-MHE-extrac-condition}
holds true if the involved $\K$ functions $\{\mu_{1},\,\mu_{2},\,\mu_{3},\,\alpha_{1},\,\alpha_{2}\}$,
$\Kinf$ functions $\{\underbar{\ensuremath{\gamma}}_{w},\thinspace\underbar{\ensuremath{\gamma}}_{v}\}$
and $\LL$ functions $\{\varphi_{1},\,\varphi_{2}\}$ satisfy the following inequality:
\begin{multline}
\mu_{1}\left(3s+3\mu_{2}^{-1}\left(3\mu_{3}(s)\right)\right)\varphi_{1}(T_{\eta,\bar{s}})\\
+\alpha_{1}\left(3\underbar{\ensuremath{\gamma}}_{w}^{-1}\left(3\mu_{3}(s)\varphi_{2}(T_{\eta,\bar{s}})\right)\right)\\
+\alpha_{2}\left(3\underbar{\ensuremath{\gamma}}_{v}^{-1}\left(3\mu_{3}(s)\varphi_{2}(T_{\eta,\bar{s}})\right)\right)\le\eta s,\label{eq: key-stability-condition-3}
\end{multline}
for all $s\in[0,\thinspace\bar{s}]$ and $t\in\mathbb{I}_{\ge0}$. 
\end{lem}
\begin{IEEEproof}
It is basically to show that inequality \eqref{eq: key-stability-condition-3}
implies inequality \eqref{eq: RGAS-of-MHE-extrac-condition} under
the conditions of Lemma \ref{lem:The-FIE-defined}. By Lemma \ref{lem:The-FIE-defined}
and its proof, we have $\rho_{x}(s,\thinspace t)=\mu_{3}(s)\varphi_{2}(t)$
and $\bar{\beta}_{x}(s,t)=\beta\left(3s+3\mu_{2}^{-1}\left(3\mu_{3}(s)\right),\,t\right)=\mu_{1}\left(3s+3\mu_{2}^{-1}\left(3\mu_{3}(s)\right)\right)\varphi_{1}(t)$.
Substituting these specific $\KL$ functions into \eqref{eq: key-stability-condition-3},
yields \eqref{eq: RGAS-of-MHE-extrac-condition} and hence completes
the proof. 
\end{IEEEproof}
\begin{lem}
\label{lem: RGAS MHE special case 1}Given conditions a)-c) of Lemma
\ref{lem:The-FIE-defined2}, the condition \eqref{eq: RGAS-of-MHE-extrac-condition}
holds true if in the given conditions the $\Kinf$ functions $\{\underbar{\ensuremath{\gamma}}_{w},\,\underbar{\ensuremath{\gamma}}_{v}\}$
and the $\K$ functions $\{\alpha_{1},\thinspace\alpha_{2}\}$ satisfy
the following inequality: 
\begin{multline}
c_{1}\left[3\left(1+^{a_{2}}\sqrt{3}\right)\right]{}^{a_{1}}s^{a_{1}}\left(T_{\eta,\bar{s}}+1\right)^{-b_{1}}\\
+\alpha_{1}\left(3\underbar{\ensuremath{\gamma}}_{w}^{-1}\left(3c_{2}s^{a_{2}}\left(T_{\eta,\bar{s}}+1\right)^{-b_{2}}\right)\right)\\
+\alpha_{2}\left(3\underbar{\ensuremath{\gamma}}_{v}^{-1}\left(3c_{2}s^{a_{2}}(T_{\eta,\bar{s}}+1)^{-b_{2}}\right)\right)\le\eta s,\label{eq: key-stability-condition-3-1}
\end{multline}
for all $s\in[0,\thinspace\bar{s}]$ and $t\in\mathbb{I}_{\ge0}$. 
\end{lem}
\begin{IEEEproof}
With $\mu_{1}(s)=c_{1}s^{a_{1}}$, $\mu_{2}(s)=\mu_{3}(s)=c_{2}s^{a_{2}}$, $\varphi_{1}(t)=(t+1)^{-b_{1}}$
and $\varphi_{2}(t)=(t+1)^{-b_{2}}$, it is straightforward to show
that \eqref{eq: key-stability-condition-3-1} is equivalent to \eqref{eq: key-stability-condition-3}.
The conclusion follows immediately from Lemma \ref{lem: RGAS of MHE - KdL form}. 
\end{IEEEproof}
Note that, to satisfy inequality \eqref{eq: key-stability-condition-3-1},
the parameter $a_{1}$ of the i-IOSS system must satisfy $a_{1}\ge1$. 
\begin{lem}
\label{lem: RGAS MHE special case 2}Given conditions a)-c) of Lemma
\ref{lem:The-FIE-defined3}, the condition \eqref{eq: RGAS-of-MHE-extrac-condition}
holds true if in the given conditions the $\Kinf$ functions $\{\underbar{\ensuremath{\gamma}}_{w},\,\underbar{\ensuremath{\gamma}}_{v}\}$
and the $\K$ functions $\{\alpha_{1},\thinspace\alpha_{2}\}$ satisfy
the following inequality: 
\begin{multline}
c_{1}\left[3\left(1+^{a_{2}}\sqrt{3}\right)\right]{}^{a_{1}}s^{a_{1}}b_{1}^{T_{\eta,\bar{s}}}+\alpha_{1}\left(3\underbar{\ensuremath{\gamma}}_{w}^{-1}\left(3c_{2}s^{a_{2}}b_{2}^{T_{\eta,\bar{s}}}\right)\right)\\
+\alpha_{2}\left(3\underbar{\ensuremath{\gamma}}_{v}^{-1}\left(3c_{2}s^{a_{2}}b_{2}^{T_{\eta,\bar{s}}}\right)\right)\le\eta s,\label{eq: key-stability-condition-3-2}
\end{multline}
for all $s\in[0,\thinspace\bar{s}]$ and $t\in\mathbb{I}_{\ge0}$. 
\end{lem}
\begin{IEEEproof}
With $\mu_{1}(s)=c_{1}s^{a_{1}}$, $\mu_{2}(s)=\mu_{3}(s)=c_{2}s^{a_{2}}$, $\varphi_{1}(t)=b_{1}^{t}$
and $\varphi_{2}(t)=b_{2}^{t}$, it is straightforward to show that
\eqref{eq: key-stability-condition-3-2} is equivalent to \eqref{eq: key-stability-condition-3}.
The conclusion follows immediately from Lemma \ref{lem: RGAS of MHE - KdL form}. 
\end{IEEEproof}

{
Inequalities \eqref{eq: key-stability-condition-3}-\eqref{eq: key-stability-condition-3-2} are materialization of inequality \eqref{eq: RGAS-of-MHE-extrac-condition} under more specific conditions on the i-IOSS property of the system and the cost function of the MHE. Therefore, the remark and interpretation on \eqref{eq: RGAS-of-MHE-extrac-condition} which are given after the proof of Theorem \ref{thm: RGAS-of-MHE} are applicable to these three inequalities.
}

Analog to the case of FIE, a specific sub-cost function $V_{T,2}$
for MHE to satisfy the conditions of Lemma \ref{lem: RGAS MHE special case 1}
or \ref{lem: RGAS MHE special case 2} can take the following form:
\begin{align}
 & V_{T,2}(\boldsymbol{\omega}_{t-T:t-1},\,\boldsymbol{\nu}_{t-T:t}) \nonumber\\
 & :=\dfrac{\lambda_{w}}{T}\sum_{\tau\in\mathbb{I}_{t-T:t-1}}l_{w,\tau}(\omega_{\tau})+\dfrac{\lambda_{v}}{T+1}\sum_{\tau\in\mathbb{I}_{t-T:t}}l_{v,\tau}(\nu_{\tau}) \nonumber\\
 & +(1-\lambda_{w})\max_{\tau\in\mathbb{I}_{t-T:t-1}}l_{w,\tau}(\omega_{\tau})+(1-\lambda_{v})\max_{\tau\in\mathbb{I}_{t-T:t}}l_{w,\tau}(\nu_{\tau}), \label{eq: specific-form-of-MHE-cost-function}
\end{align}
with given constants $\lambda_{w},\,\lambda_{v}\in[0,\,1)$. Here
the functions $l_{w,\tau}$ and $l_{v,\tau}$ are bounded as per \eqref{eq: specific-subcost-form - b},
and the resulting $\Kinf$ bound functions $\{\underbar{\ensuremath{\gamma}}_{w},\,\underbar{\ensuremath{\gamma}}_{v}\}$
which are associated with $V_{T,2}$ satisfy inequality \eqref{eq: key-stability-condition-3-1}
(or \eqref{eq: key-stability-condition-3-2}) of Lemma \ref{lem: RGAS MHE special case 1}
(or \ref{lem: RGAS MHE special case 2}). 

\begin{rem} \label{rem: on-max-term}
{
It can be shown that, if the sub-cost $V_{T,1}$ admits a form which decays with a higher order of $T$ than $V_{T,2}$ in \eqref{eq: specific-form-of-MHE-cost-function} does, then  the MHE remains RGAS even if the weight parameters $\lambda_{w}$ and $\lambda_{v}$ take the value of 1 (i.e., no max terms exist in $V_{T,2}$ in \eqref{eq: specific-form-of-MHE-cost-function}). However, the $\K$ functions in the resulting RGAS property will be dependent on the size of the moving horizon ($T$) implemented in MHE. Motivated by a relevant proof in \cite{muller2017nonlinear}, the proof can be developed by showing that $|\hat{w}_\tau|$, $\forall \tau \in \mathbb{I}_{t-T:t-1}$ and consequently $\left\Vert \hat{\boldsymbol{w}}_{t-T:t-1}\right\Vert$ is upper bounded by a sum of $\K$ functions that are dependent on $T$, which is likewise applicable to $|\hat{v}_\tau|$, $\forall \tau \in \mathbb{I}_{t-T:t-1}$. The remaining proof is first to prove the RLAS of MHE (cf. the beginning of Section \ref{sec: RGAS of the MHE}) by following the routine of the proof for the RGAS of FIE (refer to \cite{hu2015optimization}), and then use the result to establish the RGAS of MHE by following the routine of the proof of Lemma \ref{lem: MHE-FIE}. The conclusion can be generalized by assuming $V_{T,1}$ and $V_{T,2}$ to be general $\KL$ functions. Similar conclusions, however, are not proved for the associated FIE.
}
\end{rem}

\section{Numerical Examples\label{sec:numerical example}}

This section applies MHE to estimate the states of a linear system
and a nonlinear system. The two systems are provable to be i-IOSS,
and were subject to Gaussian disturbances, each of which was truncated to
the range of $[-3\sigma,\thinspace3\sigma]$, with $\sigma^{2}$ representing
the variance of the disturbance. 
{
In MHE, the prior $\bar{x}_{t-T}$ of state $x_{t-T}$ is chosen to be equal to the past MHE estimate $\hat{x}_{t-T}^{\star}$ for all $t\ge T+1$, which makes Assumption \ref{assump: A3} always satisfied. 
}%
{
The optimization problems in MHE were solved in MATLAB (version R2010b) which ran on a laptop with Intel(R) Core(TM) i7-6700HQ and CPU@2.60 GHz. Specifically, in both examples the optimization problems were solved by the ``fmincon'' solver which implements an interior-point algorithm. The iterations were set large enough such that the optimal estimates were returned.
}

The performances of MHE are compared
with those of KF (with dynamic gains) in the linear case and EKF in the nonlinear case. The performance was evaluated by mean error, and mean absolute error
(MAE) of the estimation which is defined below: 
\begin{align*}
{\rm MAE} & =\frac{1}{N(t_{f}+1)}\sum_{i=1}^{N}\sum_{t=0}^{t_{f}}\sum_{j=1}^{n}|x_{t,j}^{(i)}-\hat{x}_{t,j}^{(i)}|,
\end{align*}
where $t_{f}$ is the simulation duration, $N$ is the number of random instances
of the initial state and the disturbance sequence, $x_{t,j}^{(i)}$ is the $j$th state
of a state vector $x_{t}^{(i)}$ for time $t$ in instance $i$, and
$\hat{x}_{t,j}^{(i)}$ denotes the corresponding estimate. In both
examples, we set $N=100$ and $t_{f}=60$.

\subsection{A linear system}

Consider a linear discrete-time system described by:
\begin{align*}
&\left[\begin{array}{c}
x_{1,t+1}\\
x_{2,t+1}\\
x_{3,t+1}
\end{array}\right] =\left[\begin{array}{ccc}
0.74 & 0.21 & -0.25\\
0.09 & 0.86 & -0.19\\
-0.09 & 0.18 & 0.50
\end{array}\right]\left[\begin{array}{c}
x_{1,t}\\
x_{2,t}\\
x_{3,t}
\end{array}\right]\\
& \quad+\left[\begin{array}{c}
w_{1,t}\\
w_{2,t}\\
w_{3,t}
\end{array}\right],
\quad y_{t} =0.1x_{1,t}+2x_{2,t}+x_{3,t}+v_{t},\thinspace\thinspace\forall t\ge0.
\end{align*}
The disturbances $\{w_{1,t},\thinspace\thinspace w_{2,t},\thinspace\thinspace w_{3,t}\}$
and noise $\{v_{t}\}$ are four sequences of independent, zero mean,
truncated Gaussian noises with variances given by $\sigma_{w_{1}}^{2}=\sigma_{w_{2}}^{2}=\sigma_{w_{3}}^{2}=\sigma_{w}^{2}=0.04$
and $\sigma_{v}^{2}=0.01$, respectively. The initial state $x_{0}$ is a random variable independent
of the disturbances and noise, and follows a Gaussian distribution
with a mean of $\bar{x}_{0}$ and the variances of the three elements
are all given by $\sigma_{0}^{2}:=1$. The prior estimate of the
initial state is given as $\bar{x}_{0}=[1\thinspace\thinspace1\thinspace\thinspace-1]^{\top}$.

The system is exp-i-IOSS by Lemma \ref{lem: LTI-variant exp-i-IOSS}
established in the Appendix. By the definition of an i-IOSS system
in \eqref{eq:definition - i-IOSS}, the $\K\cdot\LL$ and $\K$ bound
functions are obtained as $\beta(s,\thinspace t)=3.04s\cdot0.9{}^{t}$,
$\alpha_{1}(s)=30.3s$ and $\alpha_{2}(s)\equiv0$. By Lemma \ref{lem: RGAS MHE special case 2}
and Theorem \ref{thm: RGAS-of-MHE}, for the MHE to be RGAS we can
specify its cost function as{ 
\begin{align}
 & V_{T}(\chi_{t-T},\thinspace\boldsymbol{\omega}_{t-T:t-1},\thinspace\boldsymbol{\nu}_{t-T:t})\nonumber \\
 & :=\frac{|\chi_{t-T}-\hat{x}_{t-T}|^{2}b_{2}^{T}}{\sigma_{0}^{2}}+V_{T,2}(\boldsymbol{\omega}_{t-T:t-1},\thinspace\boldsymbol{\nu}_{t-T:t})\label{eq: LE - Vt}
\end{align}
with
 \begin{align}
 & V_{T,2}(\boldsymbol{\omega}_{t-T:t-1},\thinspace\boldsymbol{\nu}_{t-T:t})\nonumber \\
 & :=\dfrac{\lambda_{w}}{\sigma_{w}^{2}T}\sum_{\tau=t-T}^{t-1}|\omega_{\tau}|^{2}+\dfrac{\lambda_{v}}{\sigma_{v}^{2}(T+1)}\sum_{\tau=t-T}^{t}|\nu_{\tau}|^{2}\nonumber \\
 & +\dfrac{1-\lambda_{w}}{\sigma_{w}^{2}}\|\boldsymbol{\omega}_{t-T:t-1}\|^{2}+\dfrac{1-\lambda_{v}}{\sigma_{v}^{2}}\|\boldsymbol{\nu}_{t-T:t}\|^{2}.\label{eq: LE - Vt-bar}
\end{align}
}The resulting MHE is named as MHE I, to distinguish it from another
MHE defined later. With this choice of cost function, the $\Kinf$
bound functions associated with the disturbances are derived as $\underbar{\ensuremath{\gamma}}_{w}(s)=(1-\lambda_{w})s^{2}/\sigma_{w}^{2}$,
$\gamma_{w}(s)=s^{2}/\sigma_{w}^{2}$, $\underbar{\ensuremath{\gamma}}_{v}(s)=(1-\lambda_{v})s^{2}/\sigma_{v}^{2}$,
and $\ensuremath{\gamma}_{v}(s)=s^{2}/\sigma_{v}^{2}$. To satisfy
the conditions of Lemma \ref{lem: RGAS MHE special case 2}, it suffices
to choose $b_{2}=0.81$ and $\lambda_{w},\,\lambda_{v}\in[0,\,1)$,
satisfying $0.9^{2}\le b_{2}<1$. Given a moving horizon size specified
by $T$, we solve the MHE subject to $\left\Vert \chi_{0}\right\Vert \le3\sigma_{0}$,
$\left\Vert \boldsymbol{\omega}_{t-T:t-1}\right\Vert \le3\sigma_{w}$
and $\left\Vert \boldsymbol{\nu}_{t-T:t}\right\Vert \le3\sigma_{v}$,
yielding the state estimate for each $t\in\mathbb{I}_{0:t_{f}}$.

\begin{figure}
\begin{centering}
\includegraphics[scale=0.43]{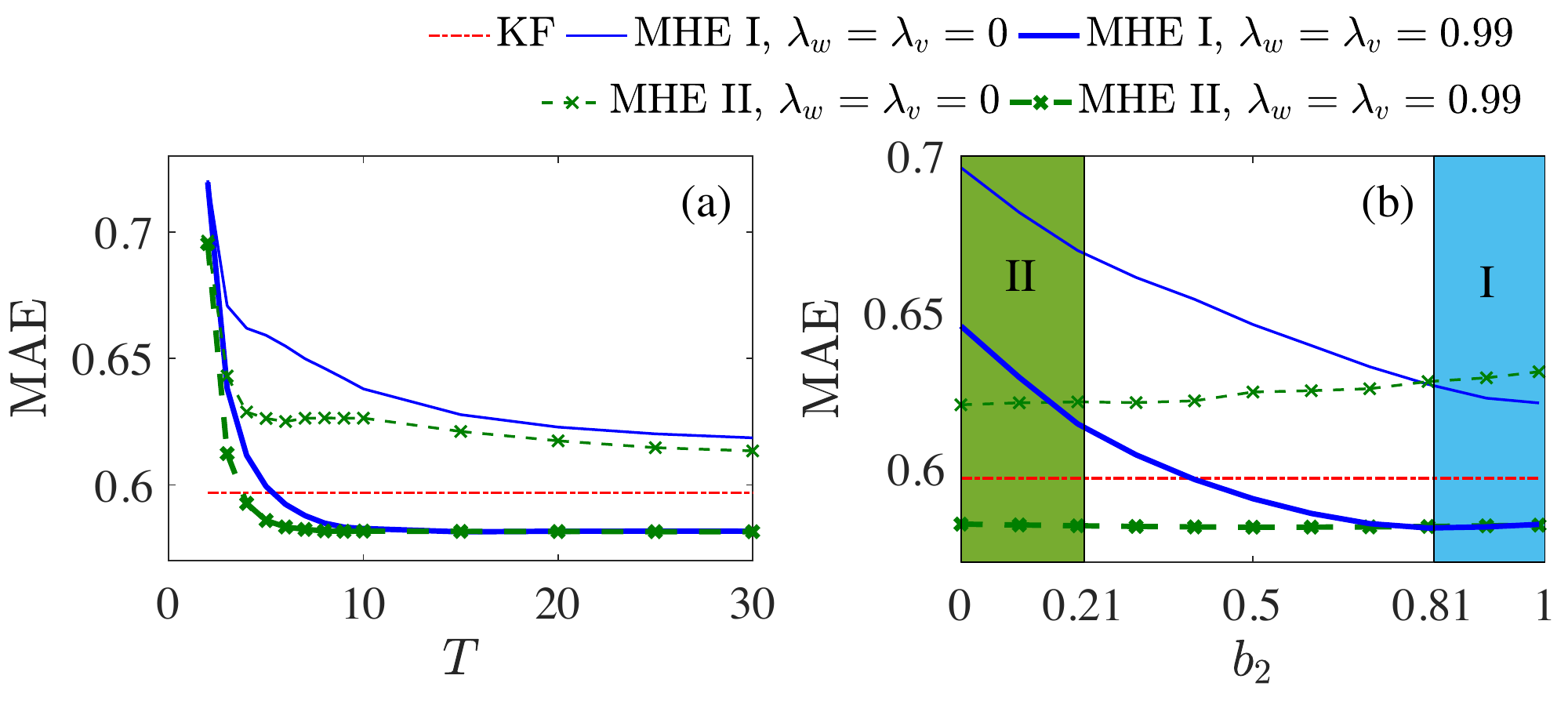} 
\par\end{centering}
\caption{MAE performances of MHE for different values of $T$ and $b_{2}$. In (a), MHE
I was implemented with $b_{2}=0.81$ and MHE II with $b_{2}=0.21$,
and in (b) both MHEs had $T=15$. The shaded areas indicate the
ranges of $b_{2}$ that meet the stability criteria of MHEs I and
II as given in Lemmas \ref{lem: RGAS MHE special case 1}-\ref{lem: RGAS MHE special case 2}.}
\label{fig: LE - MAE-T-b2 } 
\end{figure}

\begin{figure}
\begin{centering}
\includegraphics[scale=0.425]{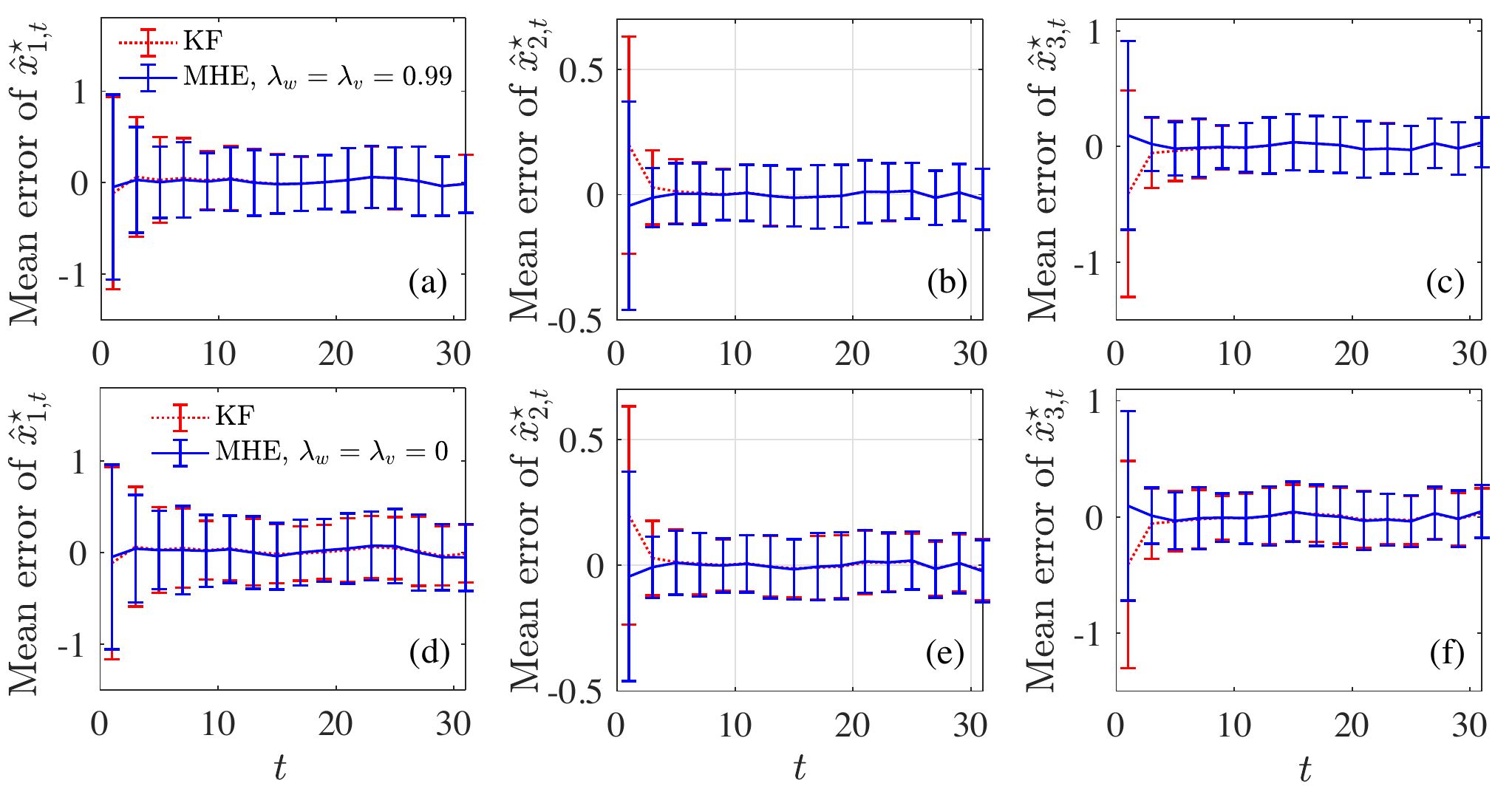} 
\par\end{centering}
\caption{Mean error performances of KF and MHE I. MHE I was implemented with
$\lambda_{w}$ and $\lambda_{v}$ both equal to 0.99 and 0 in (a)-(c) and in (d)-(f), respectively. In both cases, the two parameters
were set as $T=15$ and $b_{2}=0.81$. The bars indicate the variations
bounded by the standard deviations.}
\label{fig: LE - exemplary ME-t} 
\end{figure}

The MAEs of the estimates when $b_{2}$ and $T$ took different values
are shown in Fig. \ref{fig: LE - MAE-T-b2 }. As observed, MHE I with $\lambda_{w}=\lambda_{v}=0.99$ outperformed KF if the
horizon size $T$ and the parameter $b_{2}$ were large enough. Whereas,
MHE I was inferior to KF when $\lambda_{w}=\lambda_{v}=0$, regardless
of the values of $T$ and $b_{2}$. The observations are verified
by the results of a random instance, as shown in Fig. \ref{fig: LE - exemplary ME-t}.
We see that MHE I outperformed KF during the early stage of estimation
and became almost equivalent to KF afterwards. In addition, the results
in Fig. \ref{fig: LE - MAE-T-b2 } indicate that a small moving horizon
size, e.g., $T=10$, is sufficient for the MHE to offer a competitive
estimation, and that the improvement in the estimation performance
is marginal once the horizon is large enough. The feasible size can
thus be smaller than the sufficient size predicted by Lemma \ref{lem: RGAS MHE special case 2},
which is 39 for $\lambda_{w}=\lambda_{v}=0$ and 57 for $\lambda_{w}=\lambda_{v}=0.99$.

Alternatively, if the i-IOSS property is expressed by using a looser
$\K\cdot\LL$ bound with $\beta(s,t):=3.04s\cdot(t+1)^{\ln0.9}$ (cf.
Lemma \ref{lem: LTI-variant exp-i-IOSS} and the remark that follows),
then by Lemma \ref{lem: RGAS MHE special case 1} a different valid
cost function can be defined as: 
\begin{multline}
V_{T}'(\chi_{t-T},\thinspace\boldsymbol{\omega}_{t-T:t-1},\thinspace\boldsymbol{\nu}_{t-T:t})=\frac{|\chi_{t-T}-\hat{x}_{t-T}|^{2}}{\sigma_{0}^{2}(T+1)^{b_{2}}}\\
+V_{T,2}(\boldsymbol{\omega}_{t-T:t-1},\thinspace\boldsymbol{\nu}_{t-T:t}),\label{eq: LE - alternative Vt}
\end{multline}
where $V_{T,2}(\boldsymbol{\omega}_{t-T:t-1},\thinspace\boldsymbol{\nu}_{t-T:t})$
is the same as in \eqref{eq: LE - Vt-bar}. To satisfy the stability
conditions in Lemma \ref{lem: RGAS MHE special case 1}, it is sufficient
to choose $b_{2}=0.21$, which satisfies $0<b_{2}\le-2\ln0.9$. We
call the resulting MHE as MHE II. Simulations were performed on the
same random instances for different values of $T$ and $b_{2}$, and
the results are again shown in Fig. \ref{fig: LE - MAE-T-b2 }. The
state estimation results are slightly better than those obtained by MHE
I for different values of $T$. Similar observations were yielded
when the parameters $b_{2}$'s of the two MHEs took values in the ranges identified by Lemmas \ref{lem: RGAS MHE special case 1}
and \ref{lem: RGAS MHE special case 2}. Simulations also showed that
both MHE I and MHE II remained stable even if $b_{2}$ took values
beyond the identified ranges, which indicates the sufficiency but
non-necessity of the derived stability conditions.

The solved optimizations are convex in both MHEs. The solution times averaged over the whole simulation period (i.e., 60 time units) and 100 random instances are summarized in Table \ref{tbl: computational-times-example1}, for different sizes of moving horizons. The average solution times were less than 1.4 secs for both MHE I and MHE II if the parameters $\lambda_w$ and $\lambda_v$ were set to 0.99, and increased if $\lambda_w$ and $\lambda_v$ were set to 0 as the optimization became more challenging to solve. Moreover, in each case the solution time increased with the size of the moving horizon.

\begin{table}
\caption{Average solution time (in secs) for the linear system.}  \label{tbl: computational-times-example1}

\centering{}%
\begin{tabular}{p{2.7cm}|llllll}
\hline 
Moving horizon size & 5 & 10 & 15 & 20 & 25 & 30\tabularnewline
\noalign{\vskip3pt}
\hline 
Time for MHE I, with\newline$\lambda_{w}=\lambda_{v}=0.99$ & 0.12 & 0.30 & 0.49 & 0.72 & 1.00 & 1.33\tabularnewline\noalign{\vskip4pt}

Time for MHE I, with\newline$\lambda_{w}=\lambda_{v}=0$ & 0.60 & 1.48 & 2.24 & 2.93 & 3.49 & 3.92\tabularnewline\noalign{\vskip4pt}

Time for MHE II, with\newline$\lambda_{w}=\lambda_{v}=0.99$ & 0.10 & 0.30 & 0.45 & 0.66 & 0.91 & 1.19\tabularnewline\noalign{\vskip4pt}

Time for MHE II, with\newline$\lambda_{w}=\lambda_{v}=0$ & 0.60 & 1.46 & 2.11 & 2.58 & 3.16 & 3.79\tabularnewline
\hline 
\end{tabular}
\end{table}

Next, we compare the performances of MHE and KF when the measurements
had outliers. In this case, the noise was a mixture of two truncated
Gaussian noises: a nominal noise had a variance of $\sigma_{v}^{2}$
which occurred with a probability of $p$, and an intermittent large
noise had a variance of $100\sigma_{v}^{2}$ which occurred with a
probability of $(1-p)$ \cite{aravkin2016generalized}. The system
disturbances $w_{1,t}$, $w_{2,t}$ and $w_{3,t}$ were generated as the same
Gaussian disturbance with a variance of $\sigma_{w}^{2}$. In the
simulations, we set $p=0.9$, $\sigma_{v}=0.1$ and $\sigma_{w}=0.02$.
The other simulation settings were the same as before. In this case,
the cost function of the MHE was specified as{ 
\begin{multline*}
V_{T}(\chi_{t-T},\thinspace\boldsymbol{\omega}_{t-T:t-1},\thinspace\boldsymbol{\nu}_{t-T:t}):=\frac{|\chi_{t-T}-\hat{x}_{t-T}|^{2}b_{2}^{T}}{\sigma_{0}^{2}}\\
+\dfrac{\lambda_{w}}{\sigma_{w}^{2}T}\sum_{\tau=t-T}^{t-1}|\omega_{\tau}|^{2}+\dfrac{\lambda_{v}}{\sigma_{v}(T+1)}\sum_{\tau=t-T}^{t}|\nu_{\tau}|\\
+\dfrac{1-\lambda_{w}}{\sigma_{w}^{2}}\|\boldsymbol{\omega}_{t-T:t-1}\|^{2}+\dfrac{1-\lambda_{v}}{\sigma_{v}}\max_{\tau\in[t-T,\thinspace t]}|\nu_{\tau}|,
\end{multline*}
}which imposes 1-norm instead of 2-norm penalties on the fitting
errors $\{\nu_{\tau}\}_{t-T\le\tau\le t}$ in order to account for
outliers in the measurements in a better manner \cite{ke2005robust,hedengren2015overview,aravkin2016generalized}.
The MHE also incorporates the knowledge of identical disturbances
by including the equality constraints $\omega_{1,\tau}=\omega_{2,\tau}=\omega_{3,\tau}$
for all $\tau\in\mathbb{I}_{t-T:t-1}$. In contrast, the KF fully
implemented 2-norm penalties, and the knowledge of identical disturbances
was incorporated by specifying their covariance matrix as $\sigma_{w}^{2}\boldsymbol{1}_{3}$,
where $\boldsymbol{1}_{3}$ is a $3\times3$ matrix with all elements
equal to 1.

\begin{figure}
\begin{centering}
\includegraphics[scale=0.425]{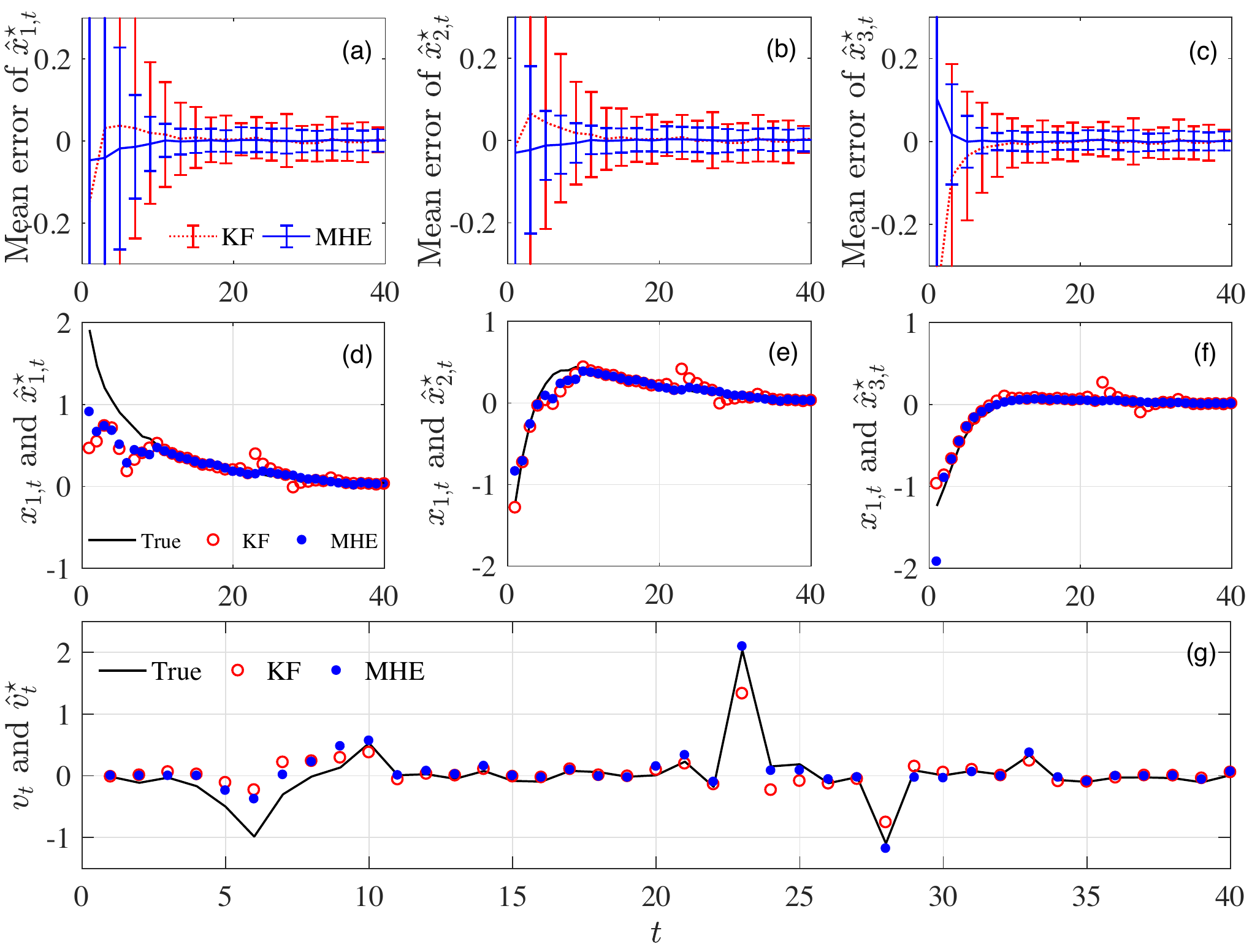} 
\par\end{centering}
\caption{Performances of MHE and KF based on measurements with outliers. Subplots (a)-(c) present the mean estimation errors, and
(d)-(g) show the estimates for a random instance.}
\label{fig: LE - mixed gaussian noise} 
\end{figure}

As shown in Fig. \ref{fig: LE - mixed gaussian noise}(a)-(c), the
state estimates obtained by MHE outperformed those obtained by KF
during the whole simulation period, in terms of both mean and variance
of the estimation errors. Fig. \ref{fig: LE - mixed gaussian noise}(d)-(f)
show the true trajectories of the two states and their associated
KF and MHE estimates in a random instance. The results confirm the
superiority of MHE in this case. Indeed, this owes to the more accurate
recovery of the measurement noise sequence by means of the 1-norm
penalties applied on the measurement fitting errors, which is supported
by the noise estimates shown in Fig. \ref{fig: LE - mixed gaussian noise}(g).

\subsection{A nonlinear system} \label{subsec: nonlinear example}

Consider a nonlinear continuous-time system described by 
\begin{equation}
\begin{aligned}\left[\begin{array}{c}
\dot{x}_{1,t}\\
\dot{x}_{2,t}
\end{array}\right] & =\left[\begin{array}{c}
-2kx_{1,t}^{2}+w_{t}\\
kx_{1,t}^{2}
\end{array}\right],\\
y_{t} & =x_{1,t}+x_{2,t}+v_{t},\thinspace\thinspace\forall t\ge0,
\end{aligned}
\label{eq: nonlinear-system}
\end{equation}
where $k=0.16$. When $w_{t}$ is constantly zero, the system describes
an ideal gas-phase irreversible reaction in a well mixed, constant
volume, isothermal batch reactor, where $x_{1,t}$ and $x_{2,t}$
represent the partial pressures and $y_{t}$ the reactor pressure
measurement \cite{haseltine2005critical,ji2016robust}. In normal
operations, the states and measurements are non-negative, i.e., $x_{1,t},\thinspace x_{2,t},\thinspace y_{t}\ge0$
for all $t\ge0$. We assume that $x_{2,0}\ge c_{0}$ for a certain
positive constant $c_{0}$. This implies that $x_{2,t}\ge c_{0}>0$
for all $t\ge0$ because $x_{2,t}$ increases with $t$.

First we prove that the system is i-IOSS, which was often assumed
without a proof in the literature, e.g., \cite{ji2016robust}. Given
two initial conditions $x_{0}^{(1)}:=[x_{1,0}^{(1)}\thinspace\thinspace x_{2,0}^{(1)}]^{\top}$
and $x_{0}^{(2)}:=[x_{1,0}^{(2)}\thinspace\thinspace x_{2,0}^{(2)}]^{\top}$,
let the corresponding state trajectory be denoted as $x_{t}^{(1)}$
and $x_{t}^{(2)}$. Define $\delta x_{1,t}=x_{1,t}^{(1)}-x_{1,t}^{(2)}$
and $p_{t}=|\delta x_{1,t}|$. The dynamics of $p_{t}$ is then derived
as 
\begin{align*}
\dot{p}_{t} & =\dfrac{\delta x_{1,t}}{|\delta x_{1,t}|}\delta\dot{x}_{1,t}=\dfrac{\delta x_{1,t}}{|\delta x_{1,t}|}\left(-2k(x_{1,t}^{(1)}+x_{1,t}^{(2)})\delta x_{1,t}+\delta w_{t}\right)\\
 & =-2k(x_{1,t}^{(1)}+x_{1,t}^{(2)})|\delta x_{1,t}|+\dfrac{\delta x_{1,t}}{|\delta x_{1,t}|}\delta w_{t}\\
 & \le-2kc_{0}p_{t}+|\delta w_{t}|.
\end{align*}
By the comparison lemma (Lemma 3.4 of \cite{khalil2002nonlinear}),
it follows that 
\begin{align*}
|\delta x_{1,t}| & =p_{t}\le p_{0}e^{-2kc_{0}t}+\int_{0}^{t}e^{-2kc_{0}(t-\tau)}|\delta w_{\tau}|d\tau\\
 & \le|\delta x_{1,0}|e^{-2kc_{0}t}+\dfrac{1-e^{-2kc_{0}t}}{2kc_{0}}\left\Vert \delta\boldsymbol{w}_{0:t}\right\Vert \\
 & \le|\delta x_{0}|e^{-2kc_{0}t}+\frac{\left\Vert \delta\boldsymbol{w}_{0:t}\right\Vert }{2kc_{0}}.
\end{align*}
Therefore, the gap between the two full states is bounded as follows:
\begin{align*}
 & |x_{t}^{(1)}-x_{t}^{(2)}|\le\left|\delta x_{1,t}\right|+\left|\delta x_{2,t}\right|\le2\left|\delta x_{1,t}\right|+\left|\delta(y_{t}-v_{t})\right|\\
 & \le2|\delta x_{0}|e^{-2kc_{0}t}+\dfrac{\left\Vert \delta\boldsymbol{w}_{0:t}\right\Vert }{kc_{0}}+\left\Vert \delta(\boldsymbol{y}_{0:t}-\boldsymbol{v}_{0:t})\right\Vert \\
 & =:\beta(|\delta x_{0}|,\thinspace t)+\alpha_{1}\left(\left\Vert \delta\boldsymbol{w}_{0:t}\right\Vert \right)+\alpha_{2}\left(\left\Vert \delta(\boldsymbol{y}_{0:t}-\boldsymbol{v}_{0:t})\right\Vert \right),
\end{align*}
where $\beta\in\KL$ and $\alpha_{1},\alpha_{2}\in\K$. Thus, the
system described in \eqref{eq: nonlinear-system} is i-IOSS by Definition
\ref{def: i-IOSS}.

\begin{figure}
\begin{centering}
\includegraphics[scale=0.53]{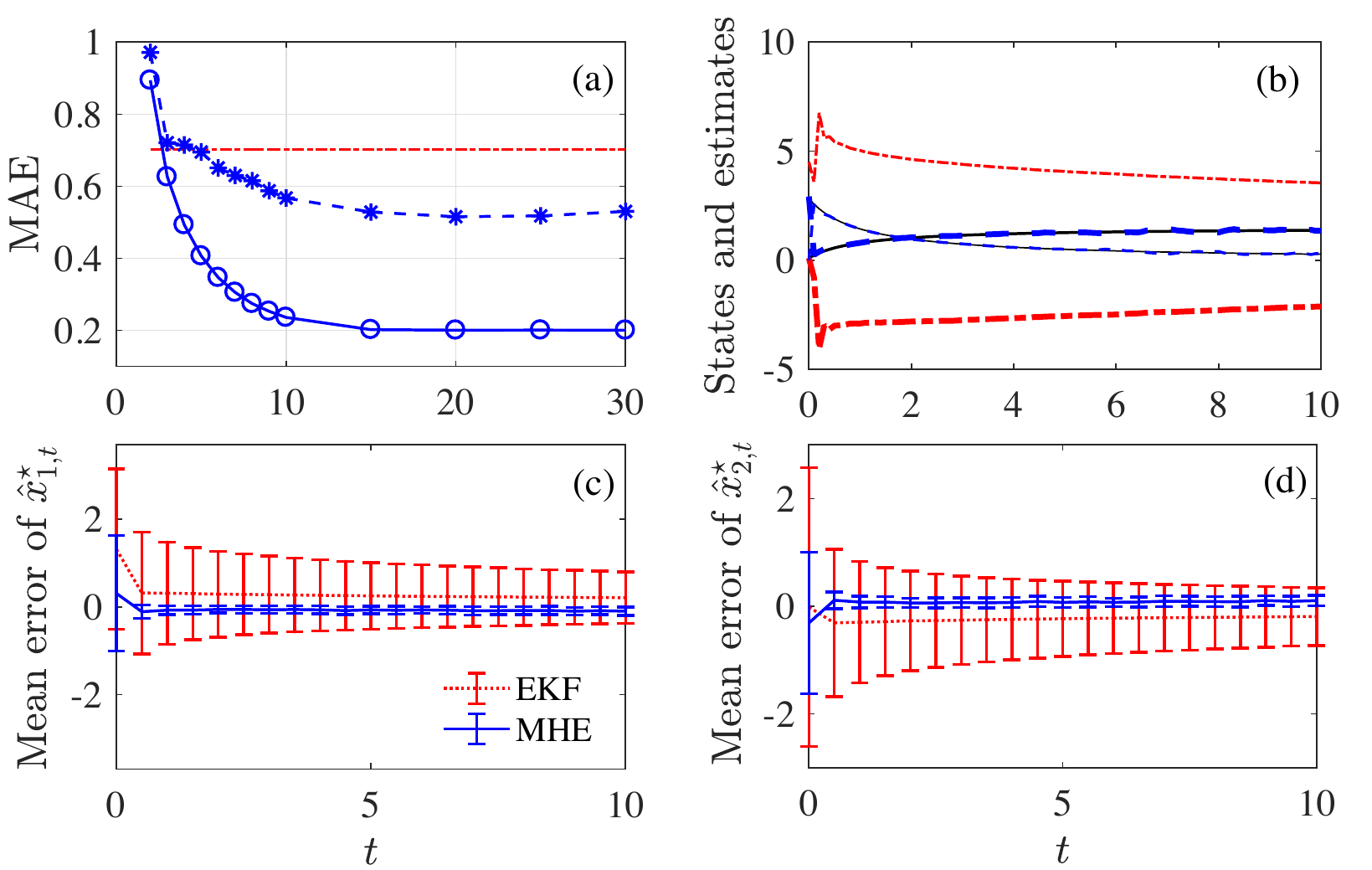} 
\par\end{centering}
\caption{MAE and mean error performances of EKF and MHE. (a) MAE performances: the red dash-dot
line is the results of EKF, and the blue dash curve with stars and the blue solid curve with circles
correspond to the results of MHE with $\lambda_{w}$ and $\lambda_{v}$
both equal to 0, and 0.99. (b)-(d) The estimation results of a random
instance: in (b) the black solid curves represent the true states, the blue
dash and red dash-dot curves correspond to the MHE and the EKF estimates, and the thinner and the thicker curves refer to states $x_{1}$ and $x_{2}$,
respectively; and mean errors of the two state estimates are shown in (c)-(d).}
\label{fig: NLE - MAE-T} 
\end{figure}

Let $w_{t}$ and $v_{t}$ be Gaussian white noises with variances
equal to $0.001^{2}$ and $0.01^{2}$, respectively. And let the initial
state $x_{0}$ follow a Gaussian distribution with a mean of $\bar{x}_{0}$
and a covariance of $\sigma_{0}^{2}I_{2}$, where $\bar{x}_{0}:=[0.1,\thinspace4.5]$,
$\sigma_{0}^{2}=9$ and $I_{2}$ is a $2\times2$ identify matrix.
{
In the simulations, we applied the Euler-Maruyama method \cite{higham2001algorithmic} to obtain discrete counterparts of the stochastic differential equations in \eqref{eq: nonlinear-system}, and the discretization step size is given by $T_{s}$. 
}
According to Lemma \ref{lem: RGAS MHE special case 2},
the MHE in discrete time is RGAS when the moving horizon size $T$
is large enough, if its cost function takes the form of \eqref{eq: LE - Vt}-\eqref{eq: LE - Vt-bar}
and is equipped with $b_{2}=e^{-4kc_{0}T_{s}}$. With $c_{0}=0.1$
and $T_{s}=0.1$ (smaller step sizes were found to yield similar results), Monte-Carlo simulations were performed for $T$
varying from 2 to 30 and the average state estimation results are
shown in Fig. \ref{fig: NLE - MAE-T}(a). As observed, MHE outperformed
EKF once the moving horizon size $T$ is larger than 2 for $\lambda_{w}=\lambda_{v}=0.99$
and 4 for $\lambda_{w}=\lambda_{v}=0$. The observations were reflected
in the results of a random instance as shown in Fig. \ref{fig: NLE - MAE-T}(b)-(d),
in which the prior estimate of the initial state was given by $\bar{x}_{0}$
and the MHE parameters were chosen as $\lambda_{w}=\lambda_{v}=0.99$
and $T=15$.

The computational times for solving the optimization problems in MHE with different sizes of moving horizons are summarized in Table \ref{tbl: computational-times-example2}. Since the optimizations are nonlinear and non-convex, the computational times were much longer than those in the previous example in which the optimizations are convex. This reflects on the challenge that persists and needs to be tackled in MHE for nonlinear systems.
\begin{table}

\caption{Average solution time (in secs) for the nonlinear system.}  \label{tbl: computational-times-example2}

\centering{}%
\begin{tabular}{p{2.5cm}|llllll}
\hline 
Moving horizon size & 5 & 10 & 15 & 20 & 25 & 30\tabularnewline
\noalign{\vskip3pt}
\hline 
Time for MHE with\newline$\lambda_{w}=\lambda_{v}=0.99$ & 0.32 & 0.70 & 1.21 & 1.78 & 2.62 & 3.59\tabularnewline\noalign{\vskip4pt}

Time for MHE with\newline$\lambda_{w}=\lambda_{v}=0$ & 0.90 & 2.89 & 5.45 & 8.11 & 11.57 & 14.50 \tabularnewline
\hline 
\end{tabular}
\end{table}

\section{Discussion\label{sec:Discussion}}
This section provides a brief discussion on solving the optimization problem defined in \eqref{eq: MHE} for MHE. If both the state and the measurement equations are linear, and if the bound sets are convex, then the optimization defined in \eqref{eq: MHE} is convex when a convex cost function is used. In this case, the optimization problem can be solved efficiently to global optimality using state-of-the-art convex solvers \cite{Boyd2004}, even if the MHE implements a large moving horizon size.

In practice, however, the state or the measurement equation is often nonlinear. This makes the optimization defined in \eqref{eq: MHE} non-convex, and the computation for a global optimal solution becomes time-consuming. This tends to void the application of MHE in cases where computational time is a key concern. To tackle the challenge, researchers have proposed solving \eqref{eq: MHE} for suboptimal solutions. For instance, in \cite{alessandri2008moving} the authors assume that the values of the cost function subjected to suboptimal solutions are within a fixed gap to the globally optimal costs. Yet it is unclear if there exist optimization solvers that can keep the assumption valid without violating the tight requirement on computational efficiency. In general, it remains an open challenge to ensure the RGAS of MHE when only suboptimal solutions are obtained for the series of optimization involved. Let the (global) optimal and the suboptimal solution of $(\chi_{t-T},\,\boldsymbol{\omega}_{t-T:t-1},\,\boldsymbol{\nu}_{t-T:t})$ be denoted as $(\hat{x}_{t-T},\,\hat{\boldsymbol{w}}_{t-T:t-1},\,\hat{\boldsymbol{v}}_{t-T:t})$ and $(\check{x}_{t-T},\,\check{\boldsymbol{w}}_{t-T:t-1},\,\check{\boldsymbol{v}}_{t-T:t})$, respectively. The following result may shed some light on the ways to tackle this challenge.

\begin{thm} \label{thm: RGAS-of-suboptimal-MHE}
Given Assumption \ref{assump: A3}
and any $\eta\in(0,\thinspace1)$, the MHE which implements a suboptimal solution for the optimization in \eqref{eq: MHE} is RGAS for all $T\ge T_{\eta,\bar{s}}$ if the following conditions are satisfied: a) the suboptimal solution yields a cost value which satisfies the following inequality,
\begin{multline}
V_t(\check{x}_{t-T},\,\check{\boldsymbol{w}}_{t-T:t-1},\,\check{\boldsymbol{v}}_{t-T:t})\\
\le \gamma\left(V_t(\hat{x}_{t-T},\,\hat{\boldsymbol{w}}_{t-T:t-1},\,\hat{\boldsymbol{v}}_{t-T:t})\right) \label{eq: suboptimal-MHE-condition}
\end{multline}
with a certain $\gamma\in \K$; and b) the cost function of the associated FIE satisfies Assumption \ref{assump: A1}, a modified Assumption \ref{assump: A2} and a modified inequality \eqref{eq: RGAS-of-MHE-extrac-condition}, in which the modifications are to replace the $\KL$ function $\rho_x$ and the  $\K$ functions $\gamma_w,\,\gamma_v$ used in Assumption \ref{assump: A2} and inequality \eqref{eq: RGAS-of-MHE-extrac-condition} with the $\KL$ function $\gamma\circ3\rho_x$ and $\K$ functions $\gamma\circ3\gamma_w$, $\gamma\circ3\gamma_v$, respectively.
\end{thm}
\begin{IEEEproof}
By condition a) and Assumption \ref{assump: A1}, we have 
\begin{align*}
 & \underbar{\ensuremath{\rho}}_{x}(\left|\check{x}_{0}-\bar{x}_{0}\right|,t)+\underbar{\ensuremath{\gamma}}_{w}(\left\Vert \check{\boldsymbol{w}}_{0:t-1}\right\Vert )+\underbar{\ensuremath{\gamma}}_{v}(\left\Vert \check{\boldsymbol{v}}_{0:t}\right\Vert )\\
 & \le V_{T}(\check{x}_{0}-\bar{x}_{0},\,\check{\boldsymbol{w}}_{0:t-1},\,\check{\boldsymbol{v}}_{0:t})\\
 & \le\gamma\left(V_{T}(\hat{x}_{0}-\bar{x}_{0},\,\hat{\boldsymbol{w}}_{0:t-1},\,\hat{\boldsymbol{v}}_{0:t})\right)\\
 & \le\gamma\left(V_{T}(x_{0}-\bar{x}_{0},\,\boldsymbol{w}_{0:t-1},\,\boldsymbol{v}_{0:t})\right)\\
 & \le\gamma\left(\rho_{x}(\left|x_{0}-\bar{x}_{0}\right|,\,t)+\gamma_{w}(\left\Vert \boldsymbol{w}_{0:t-1}\right\Vert )+\gamma_{v}(\left\Vert \boldsymbol{v}_{0:t}\right\Vert )\right)\\
 & \le\check{\rho}_{x}(\left|x_{0}-\bar{x}_{0}\right|,\,t)+\check{\gamma}_{w}(\left\Vert \boldsymbol{w}_{0:t-1}\right\Vert )+\check{\gamma}_{v}(\left\Vert \boldsymbol{v}_{0:t}\right\Vert )
\end{align*}
where $\check{\rho}_{x}:=\gamma\circ3\rho_{x}$, $\check{\gamma}_{w}:=\gamma\circ3\gamma_{w}$
and $\check{\gamma}_{v}:=\gamma\circ3\gamma_{v}$, which are $\KL$,
$\K$ and $\K$ functions, respectively. The lower and the eventual upper bounds of $V_{T}(\check{x}_{0}-\bar{x}_{0},\,\check{\boldsymbol{w}}_{0:t-1},\,\check{\boldsymbol{v}}_{0:t})$ are in the same forms of the counterparts obtained for a global optimal solution (refer to part (a) of the proof of Theorem 1 in \cite{hu2015optimization}). The only differences are that the left-hand side of the first inequality is expressed by the suboptimal instead of optimal solution variables, and that the right-hand side of the last inequality is described using the new functions $\check{\rho}_x$, $\check{\gamma}_w$ and $\check{\gamma}_v$, instead of $\rho_x$, $\gamma_w$ and $\gamma_v$. Consequently, under the modified Assumption \ref{assump: A2}, the proof of the associated FIE (which implements suboptimal solutions) being RGAS can be developed by following the same routine of part (a) of the proof of Theorem 1 in \cite{hu2015optimization}. The remaining proof is to show that an additional condition, as a counterpart of \eqref{eq: RGAS-of-MHE-extrac-condition}, when the MHE implements suboptimal solutions, is also satisfied. This is done by applying the modified inequality \eqref{eq: RGAS-of-MHE-extrac-condition} stated in condition b). The proof is then complete.
\end{IEEEproof}

Theorem \ref{thm: RGAS-of-suboptimal-MHE} indicates that MHE can be RGAS even if it implements suboptimal solutions for the series of optimizations that are revealed over time. To that end, each suboptimal solution needs to satisfy certain conditions, say, inequality \eqref{eq: suboptimal-MHE-condition} which requires the yielded cost value to be upper bounded by a $\K$ function of the counterpart that results from an optimal solution. Since the condition does not restrict the form of the $\K$ function, it implies flexibility in obtaining the suboptimal solutions and hence also a direction for future research.

On the other hand, we may model the system dynamics and the measurements using discrete-time linear equations in which the disturbances and noises lump all unmodeled nonlinear dynamics (including unknown external disturbances and noises). This will enable MHE to solve only convex programs, but meanwhile may sacrifice the estimation accuracy for the increased uncertainties. Future research may be conducted to design appropriate convex MHEs that balance between the estimation accuracy and the computational complexity.

In addition, as remarked after the proof of Lemma \ref{lem: MHE-FIE}, it is possible to change the size of the moving horizon of MHE online while keeping the MHE being RGAS. This will enable MHE with adaptive moving horizon which can be computationally more efficient on average as compared to the MHE that implements a moving horizon of a fixed size.

\section{Conclusion\label{sec:Conclusion}}

This paper proved the robust global asymptotic stability (RGAS) of
full information estimation (FIE) and its practical approximation,
moving horizon estimation (MHE) under general settings. The results
indicate that both FIE and MHE lead to bounded estimation errors under
mild conditions for an incrementally input/output-to-state stable
(i-IOSS) system subjected to bounded system and measurement disturbances.
The stability conditions require that the cost function to be optimized
has a property resembling the i-IOSS property of a system, but with
a higher sensitivity to the uncertainty in the initial state. The
stability of MHE additionally requires that the moving horizon is
long enough to suppress temporal propagation of the estimator errors. Under the same conditions, the MHE was also shown to converge to the true state if the disturbances converge to zero in time.

When dealing with constrained nonlinear systems, MHE has to solve
a non-convex program at each estimation point. Searching for a global
optimal solution to the program requires considerable computational
resources which are often unaffordable in applications. This problem
has motivated considerable efforts to develop robustly stable MHE
that relies on suboptimal but computationally more efficient solutions. We provided a brief discussion of this direction which may hopefully contribute to the future development of a systematic
and practical solution for this equally important problem.

\section*{Acknowledgment}
The author thanks anonymous reviewers for their valuable comments that have helped to improve the quality of the paper. The author is also grateful to Dr. Keyou You and Dr. Lihua Xie for their help during the early development of the results.

\section*{Appendix: A supporting lemma and its proof}
\begin{lem}
\label{lem: LTI-variant exp-i-IOSS} Consider a system described by
\eqref{eq:system}, where $f(x_{t},\thinspace w_{t})=Ax_{t}+g(x_{t})+w_{t}$
with $g$ being a nonlinear function, and $h$ is a linear or nonlinear
measurement function. Suppose that $A$ is diagonalizable as $P^{-1}\Lambda P$
for a certain non-singular matrix $P$ and a diagonal matrix $\Lambda$.
If the spectrum radius of $A$, denoted by $\rho(A)$, is less than
one and the nonlinear functions satisfy $|g(x_{t}^{(1)})-g(x_{t}^{(2)})|\le L|h(x_{t}^{(1)})-h(x_{t}^{(2)})|$
for all admissible $x_{t}^{(1)}$ and $x_{t}^{(2)}$ and a positive
constant $L$, then the system is exp-i-IOSS as per \eqref{eq:definition - i-IOSS},
in which the $\KL$ function can be specified as $\beta(s,\thinspace t)=|P^{-1}||P|s\rho^{t}(A)\le c_{1}s(t+1)^{-b_{1}}$
for all $s,\thinspace t\ge0$, and the $\K$ functions as $\alpha_{1}(s)=\frac{|P^{-1}||P|}{1-\rho(A)}\cdot s$
and $\alpha_{2}(s)=\frac{L|P^{-1}||P|}{1-\rho(A)}\cdot s$ for all
$s\ge0$. Here the parameters $b_{1}$ and $c_{1}$ are positive constants
satisfying $c_{1}\ge\max_{t\in\mathbb{I}_{\ge0}}|P^{-1}||P|\rho^{t}(A)(t+1)^{b_{1}}$. 
\end{lem}
\begin{IEEEproof}
Given two initial conditions $x_{0}^{(1)}$ and $x_{0}^{(2)}$
and corresponding disturbance sequences $\boldsymbol{w}_{0:t-1}^{(1)}$
and $\boldsymbol{w}_{0:t-1}^{(2)}$, let the system states at time
$t$ be yielded as $x_{t}^{(1)}$ and $x_{t}^{(2)}$, respectively.
By using the analytical expressions of the two states, we have {

\vspace{-6pt}
{\footnotesize{} 
\begin{align*}
 & |x_{t}^{(1)}-x_{t}^{(2)}|\\
= & \left|A^{t}(x_{0}^{(1)}-x_{0}^{(2)})+\sum_{l=0}^{t-1}A^{t-1-l}\left(g(x_{l}^{(1)})+w_{l}^{(1)}-g(x_{l}^{(2)})-w_{l}^{(2)}\right)\right|\\
\le & |A^{t}||x_{0}^{(1)}-x_{0}^{(2)}|+\left(\begin{array}{c}
\|\boldsymbol{w}_{0:t-1}^{(1)}-\boldsymbol{w}_{0:t-1}^{(2)}\|\\
+\|g(\boldsymbol{x}_{0:t-1}^{(1)})-g(\boldsymbol{x}_{0:t-1}^{(2)})\|
\end{array}\right)\sum_{k=0}^{t-1}|A^{k}|\\
\le & |P^{-1}||P|\left(\begin{array}{c}
\rho^{t}(A)|x_{0}^{(1)}-x_{0}^{(2)}|\\
+\left(\begin{array}{c}
\|\boldsymbol{w}_{0:t-1}^{(1)}-\boldsymbol{w}_{0:t-1}^{(2)}\|\\
+L\|h(\boldsymbol{x}_{0:t-1}^{(1)})-h(\boldsymbol{x}_{0:t-1}^{(2)})\|
\end{array}\right)\sum_{k=0}^{t-1}\rho^{k}(A)
\end{array}\right)\\
\le & |P^{-1}||P|\left(\begin{array}{c}
\rho^{t}(A)|x_{0}^{(1)}-x_{0}^{(2)}|+\dfrac{1}{1-\rho(A)}\|\boldsymbol{w}_{0:t-1}^{(1)}-\boldsymbol{w}_{0:t-1}^{(2)}\|\\
+\dfrac{L}{1-\rho(A)}\|h(\boldsymbol{x}_{0:t}^{(1)})-h(\boldsymbol{x}_{0:t}^{(2)})\|
\end{array}\right).
\end{align*}
}}The last inequality implies that the system is i-IOSS and in particular
exp-i-IOSS by Definition \ref{def: i-IOSS}, in which the $\KL$ function
$\beta$ and the $\K$ functions $\{\alpha_{1},\thinspace\alpha_{2}\}$
are specified by the lemma. Given $b_{1}\ge0$, define $c_{1,\min}^{b_{1}}=\max_{t\in\mathbb{I}_{\ge0}}|P^{-1}||P|\rho^{t}(A)(t+1)^{b_{1}}$.
Since the maximum exists while $\rho(A)<1$, $c_{1,\min}^{b_{1}}$
is well defined. Therefore, there always exists $c_{1}\ge c_{1,\min}^{b_{1}}$
such that $|P^{-1}||P|\rho^{t}(A)\le c_{1}(t+1)^{-b_{1}}$ for all
$t\ge0$. The conclusions of the lemma follows immediately.
\end{IEEEproof}{\footnotesize \par}

By the lemma, if a system is exp-i-IOSS with the $\KL$ bound function
given as $\beta(s,\thinspace t)=c_{1}s^{a_{1}}b_{1}^{t}$ for certain
positive constants $a_{1}$, $b_{1}$ and $c_{1}$ with $b_{1}<1$,
then it is always feasible to specify a looser $\KL$ bound function
as $\beta'(s,\thinspace t)=c_{1}'s^{a_{1}}(t+1)^{-b_{1}'}$ for certain
positive constants $b_{1}'$ and $c_{1}'$ satisfying $c_{1}'\ge\max_{t\in\mathbb{I}_{\ge0}}c_{1}b_{1}^{t}(t+1)^{b_{1}'}$.
For example, if $b_{1}'$ is set to $1$, then $c_{1}'$ can be specified
as $\max c_{1}b_{1}^{t}(t+1)$ which exists for $t\in\mathbb{I}_{\ge0}$
when $0\le b_{1}<1$. For another example, if $b_{1}'$ is set to
$-\ln b_{1}$, then $c_{1}'$ can be specified as equal to $c_{1}$
(which is equal to $\max_{t\in\mathbb{I}_{\ge0}}c_{1}b_{1}^{t}(t+1)^{-\ln b_{1}}$).

\bibliographystyle{IEEEtran}
\bibliography{Hu_MHE_2016}

\begin{IEEEbiography}[{\includegraphics[width=1in,height=1.25in,clip,keepaspectratio]{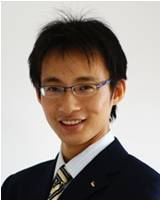}}]{Wuhua Hu} 
received the BEng degree in Automation in 2005 and the MEng degree
in Detecting Technique and Automation Device in 2007 from Tianjin
University, China. He received the PhD degree in Communication Engineering
from Nanyang Technological University (NTU), Singapore, in 2012. He
is currently a Research Scientist with the Institute for Infocomm
Research, Agency for Science, Technology and Research (A$^\star$STAR),
Singapore. Before joining A$^\star$STAR, he worked as a Research Fellow
first with the School of Mechanical and Aerospace Engineering and
then the School of Electrical and Electronic Engineering, NTU, Singapore,
from Aug 2011 to Mar 2016. His research interests lie in modeling,
estimation, control and optimization of dynamical systems, with applications
for smart and greener power and energy systems. 
\end{IEEEbiography}

\end{document}